\documentclass[11pt,a4paper]{article}

\usepackage[english]{babel}
\usepackage{amsmath,amssymb,amsfonts,a4}
\usepackage{eucal,bbm,stmaryrd}
\usepackage{hyperref}
\usepackage{xcolor,tikz}

\usepackage{apptools}
\AtAppendix{\counterwithin{Lem}{section}}

\begin{document}

\numberwithin{equation}{section}

\newtheorem{Thm}{Theorem}
\newtheorem{Prop}{Proposition}
\newtheorem{Def}{Definition}
\newtheorem{Lem}{Lemma}
\newtheorem{Rem}{Remark}
\newtheorem{Cor}{Corollary}
\newtheorem{Con}{Conjecture}

\newcommand{\Thmautorefname}{Theorem}
\newcommand{\Propautorefname}{Proposition}
\newcommand{\Defautorefname}{Definition}
\newcommand{\Lemautorefname}{Lemma}
\newcommand{\Remautorefname}{Remark}
\newcommand{\Corautorefname}{Corollary}
\newcommand{\Conautorefname}{Conjecture}

\newcommand{\Pf}{\noindent{\bf Proof: }}
\newcommand{\qed}{\hspace*{3em} \hfill{$\square$}}

\newcommand{\N}{\mathbbm{N}}
\newcommand{\Z}{\mathbbm{Z}}

\newcommand{\E}{\mathbbm{E}}
\renewcommand{\P}{\mathbbm{P}}

\newcommand{\Bin}{\mathit{Bin}}

\newcommand{\la}{\lambda}
\newcommand{\La}{\Lambda}
\newcommand{\si}{\sigma}
\newcommand{\al}{\alpha}
\newcommand{\be}{\beta}
\newcommand{\ep}{\epsilon}
\newcommand{\ga}{\gamma}
\newcommand{\Ga}{\Gamma}
\newcommand{\de}{\delta}
\newcommand{\De}{\Delta}
\newcommand{\ph}{\varphi}
\newcommand{\om}{\omega}
\renewcommand{\th}{\theta}
\newcommand{\vth}{\vartheta}
\newcommand{\ka}{\kappa}

\newcommand{\cZ}{\mathcal{Z}}
\newcommand{\cT}{\mathcal{T}}
\newcommand{\cX}{\mathcal{X}}
\newcommand{\cY}{\mathcal{Y}}
\newcommand{\cW}{\mathcal{W}}
\newcommand{\tcW}{\tilde{\cW}}
\newcommand{\bcW}{\bar{\cW}}
\newcommand{\cV}{\mathcal{V}}
\newcommand{\cC}{\mathcal{C}}
\newcommand{\cN}{\mathcal{N}}
\newcommand{\da}{\dagger}
\newcommand{\Mda}{M^{\dagger}}
\newcommand{\tM}{\tilde{M}}
\newcommand{\xda}{x_\da}
\newcommand{\xst}{x_*}

\newcommand{\Aut}{Aut}

\newcommand{\Nm}{N_{\wedge}}
\newcommand{\sq}{\times}

\newcommand{\stm}{\setminus}
\newcommand{\lra}{\leftrightarrow}
\newcommand{\lraa}[1]{\stackrel{#1}{\lra}}
\newcommand{\Ra}{\Rightarrow}
\newcommand{\Lra}{\Leftrightarrow}

\definecolor{lgray}{gray}{0.9}

\thispagestyle{plain}
\title{Monotonicity of Markov chain transition probabilities\\ via quasi-stationarity - an application to \\
Bernoulli percolation on $C_k \times \Z$}
\author{Philipp König and Thomas Richthammer
\footnote{Institut f\"ur Mathematik, Universit\"at Paderborn}}

\maketitle
\begin{abstract}
Let $X_n, n \ge 0$ be a Markov chain with finite state space $M$. 
If $x,y \in M$ such that $x$ is transient we have 
$\P^y(X_n = x) \to 0$ for $n \to \infty$, and under mild aperiodicity conditions this convergence is monotone 
in that for some $N$ we have $\forall n \ge N: 
\P^y(X_n = x)$ $\ge \P^y(X_{n+1} = x)$. 
We use bounds on the rate of convergence of the Markov chain to its quasi-stationary distribution to obtain explicit bounds on $N$. 
We then apply this result to Bernoulli percolation with parameter $p$ 
on the cylinder graph $C_k \times \Z$. 
Utilizing a Markov chain describing infection patterns layer per layer, 
we thus show the following uniform result on the monotonicity of 
connection probabilities: 
$\forall k \ge 3\,  \forall n \ge 500k^6 2^k \,\forall p \in (0,1) \,  \forall m \in C_k\!\!:$ $\P_p((0,0) \lra (m,n)) \ge 
\P_p((0,0) \lra (m,n+1))$. 
In general these kind of monotonicity properties 
of connection probabilities are difficult to establish 
and there are only few pertaining results. 
\end{abstract}

%



\section{Introduction}

Let us consider a homogeneous Markov chain $X_n, n \ge 0$, 
in discrete time with finite state space $M$. 
Suppose that we have a decomposition $M = M_0 \cup M_1$ 
such that $M_0$ is transient and $M_1$ is absorbing. 
We interpret $M_0$ as states such that the chain is alive and $M_1$ as states such that the chain is dead, and we are interested in the values of the chain only as long as it is alive. 
It is easy to see that for all $x,y \in M_0$ we have 
$\P^y(X_n = x) \to 0$ for $n \to \infty$, 
and it is natural to ask, whether this convergence is monotone
in that for some $N$ we have 
\begin{equation} \label{equ:Markovmono}
\forall n \ge N: \P^y(X_n  = x) \ge \P^y(X_{n+1} = x). 
\end{equation}
This is not too difficult to investigate, since the asymptotic behavior of the transition probabilities is well known, 
see e.g. \cite{Li2}: Under mild aperiodicity conditions 
there are constants $\ga > 0, \si \in \{1,2,...\}, 
\la \in (0,1)$ depending on $x,y$ such that 
$\P^y(X_n = x) \sim \ga n^{\si - 1} \la^n$. 
This asymptotic exponential decay readily implies 
that for some $N$ we have \eqref{equ:Markovmono}. 

From this general argument however it is not so obvious how to find an explicit value of $N$ for a given Markov chain, 
so that \eqref{equ:Markovmono} holds. 
For simplicity let us assume that the structure of the state space is as simple as possible, i.e. $M_0$ is communicating and aperiodic, 
which implies that $\P^y(X_n = x) > 0$ for sufficiently large $n$.  
Under this additional assumption 
the chain admits a unique quasi-stationary distribution $\al$, 
i.e. a distribution $\al$ supported on $M_0$ such that 
\begin{equation} \label{equ:quasi1}
\forall n \ge 0 \forall x \in M_0: \P^\al(X_n = x| X_n \in M_0) = \al(x). 
\end{equation}
$\al$ can be obtained as a left eigenvector of the 
restriction of the transition matrix of $X_n$ to $M_0$. 
The corresponding eigenvalue $\la \in (0,1)$ satisfies
$\P^\al(X_n \in M_0) = \la^n$ for all $n \ge 0$. 
The asymptotic behavior of the  transition probabilities
can now be rephrased in terms of these objects: 
For every $y \in M_0$ we have $\P^y(X_n \in M_0) \sim c_y \la^n$ for some 
constant $c_y>0$ and 
\begin{equation} \label{equ:quasi2}
\P^y(X_n = x| X_n \in M_0) \to \al(x) \qquad 
\text{ for } n \to \infty.  
\end{equation}
For the above facts see \cite{DP} 
and the references therein. For more on quasi-stationary distribution we refer to \cite{CMM}, in which however 
the focus is on Markov processes in continuous time as it is in most of the literature. 
A recent result of Champagnat and Villemonais (\cite{CV}) 
gives explicit bounds on the rate of convergence in \eqref{equ:quasi2}. 
We show how this can be used 
to get an explicit bound on $N$ such that the monotoncity condition \eqref{equ:Markovmono} holds (uniformly in $x,y \in M_0$). 

\bigskip

Our interest in the monotonicity property \eqref{equ:Markovmono}
stems from a similar monotonicity property in Bernoulli percolation. We give a brief introduction to this well known probability model: 
Let $G = (V,E)$ be a (simple) graph  and $p \in (0,1)$. 
Bernoulli percolation on $G$ with parameter $p$ is defined in 
terms of a family $\cZ_e, e \in E,$ of $\{0,1\}$-valued random variables such that 
$\cZ_e, e \in E,$ are independent and $\P_p(\cZ_e = 1) = p$. 
Here we indicate the dependency on the parameter $p$
in the notation for probabilities (and expectations). 
$e$ is called open if $\cZ_e = 1$, and $e$ is called closed if 
$\cZ_e = 0$. Objects of interest are connection events 
of the form 
$$
\{u \lra v\} := \{\exists n \ge 0, v_0,...,v_n \in V \!: 
v_0 = u, v_n = v, 
\forall i: v_iv_{i+1} \in E, \cZ_{v_iv_{i+1}} = 1\}
$$
for $u,v \in V$, and the clusters of vertices given 
by $C(v) := \{u \in V: u \lra v\}$ for $v \in V$. 
Bernoulli percolation has mostly been studied on very regular graphs such as lattices, for an overview we refer to \cite{G}. 
We will consider a fixed vertex $o \in V$, which is interpreted as the origin of an infection, 
which may spread on $G$ via open edges. 
Thus $\{o \lra v\}$ represents the event that $v$ is infected and $C(o)$ represents the sets of infected vertices. 
It is not unreasonable to expect that vertices closer to $o$ 
(in some sense that has to be specified) are more likely to be infected than vertices further away. 
In order to make this into a precise statement, we consider 
layered graphs $G^\sq$, 
which are defined to be Cartesian products of  a
graph $G = (V,E)$ and the 1D lattice $\Z$. 
Loosely speaking the layered graph $G^\sq$ consists of biinfinitely many copies of $G$ (forming the layers of $G^\sq$) 
with vertices in two distinct layers connected to each other 
iff they are copies of the same vertex in $V$ and 
the two layers are adjacent. 
To illustrate this definition we note that $\Z^d$ can be viewed 
as the layered graph $(\Z^{d-1})^\sq$. 
For Bernoulli percolation on $G^\sq$ it can be conjectured that for any fixed $o \in V$ 
\begin{equation}\label{equ:layer}
\forall v \in V \; \forall n \ge 0: 
\;
\P_p((o,0) \lra (v,n)) \ge \P_p((o,0) \lra (v,n+1)). 
\end{equation}
This kind of monotonicity property is very similar to the better known bunkbed conjecture, see \cite{Ri2} and the references therein. 
It seems that monotonicity properties of this sort are difficult to prove (if they are true at all), 
even though it is very compelling to assume that they hold. 
Indeed, there are only very few pertaining results. 
\cite{LPS} considers
Bernoulli percolation on $\Z^d = (\Z^{d-1})^\sq$
and shows that \eqref{equ:layer} holds in the special case 
$v = o$ provided that $p$ is sufficiently close to $0$. 
This result relies on the asymptotics of connection probabilities in \cite{CCC}, and it is not clear how to extract an  explicit condition on $p$ from the given proof. 
In \cite{KR} we investigate the above monotonicity 
property for general finite connected $G$ and show 
that there is a constant $N = N(G)$ such that 
\begin{equation} \label{equ:layerN}
\forall v \in V \forall n \ge N \forall p \in [0,1]\!: 
\P_p((o,0) \!\lra\! (v,n)) \ge \P_p((o,0) \!\lra\! (v,n+1)),  
\end{equation}
which is weaker than \eqref{equ:layer} in that 
we have monotonicity only from some point $N$ on (but $N$ does not depend on $v$ or $p$). This result relies on a Markov 
chain that builds the percolation pattern layer by layer and only uses open bonds below the layer under consideration. 
The main ingredient of the proof is the 
asymptotics of transition probabilities of this Markov chain. 
Again, it is not possible 
to extract an explicit condition on $n$, 
i.e. an explicit bound on $N$, from the given proof. 
Our paper will address this problem making use of quasi-stationarity as described above. 
We think it is instructive to see how to utilize the explicit bounds of \cite{CV}, how to make good choices of the control parameters introduced there, and how to make use of 
percolation specific arguments, 
in order to obtain a good bound on $N$. 
To make the result as concrete as possible we will focus on the layered graph $C_k^\sq$, 
where $C_k$ is the cycle graph on $k$ vertices. 
We note that \cite{KR} and the present paper should be considered as companion papers. 
Whereas \cite{KR} introduces the Markov chain approach 
for percolation on layered graphs and proves
for every finite graph $G$ the existence of $N$  such that \eqref{equ:layerN} holds, the focus of the present paper is on how to find such an $N$, which is as small as possible, 
for a concrete given graph $G$. 
While some of our arguments are thus similar to those of \cite{KR}, important parts of the general strategy need to be changed, 
and we need more refined combinatorial and probabilistic estimates 
in order to get a good bound on $N$. 
As in \cite{KR} we complement this result by a result on the monotonicity of the expected number of 
infected vertices per layer in case $p$ is sufficiently small, 
where again an explicit bound on $p$ is provided.

The paper is organized as follows. 
In Section 2 we present our results. 
In Section 3 we prove our result on the monotonicity of 
Markov chain transition probabilities. 
In Section 4 we state and prove some properties of the Markov chain 
of infection patterns.  
In Sections 5 we prove our result on the monotonicity 
of connection probabilities for Bernoulli percolation on 
$C_k^\sq$, and in Section 6 we prove our complementary monotonicity result on the expected number of infected vertices per layer. 
In an appendix we collect some purely analytical estimates used in the preceding chapters.

\section{Results}

We first state a general criterion for the onset of monotonicity 
for transition probabilities of Markov chains in the setting 
described in the introduction. 

\begin{Thm} \label{Thm:quasi}
Let $(X_n)_{n \ge 0}$ be a Markov chain  with 
finite state space $M_0 \cup M_1$, where $M_1$ is absorbing, 
$M_0$ is transient and communicating with $|M_0| > 1$. 
Suppose that there is a distribution $\nu$ on $M_0$ and there are constants $n_\nu \in \{1,2,...\}$ and 
$c_\nu,c_\nu' \in (0,1]$ such that 
\begin{align} \label{equ:c1}
&\forall y,x \in M_0: 
\quad 
\P^y(X_{n_\nu} = x|X_{n_\nu} \in M_0) \ge c_\nu \nu(x) \quad \text{ and } \\
\label{equ:c2}
&\forall y \in M_0, n \ge 0: \quad 
\P^\nu(X_n \in M_0) \ge c_\nu' \P^y(X_n \in M_0).
\end{align}
\begin{itemize}
\item[(a)] 
There is a unique quasi-stationary distribution $\al$ on $M_0$, 
and for every initial distribution $\mu$ on $M_0$ we have  
\begin{equation} \label{equ:exponentialquasi}
\forall n \ge 0: \sum_{x \in M_0}|\P^\mu(X_n = x|X_n \in M_0) - \al(x)| \le 2 (1-c_\nu c_\nu')^{\lfloor n/n_\nu\rfloor}. 
\end{equation}
\item[(b)]
There are constants $c_\al,c_\da \in (0,1]$ and $n_\da \in \{0,1,...\}$  such that 
\begin{align}\label{equ:cal}
&\forall x \in M_0: \al(x) \ge c_\al, \\
\label{equ:cda}
&\forall y \in M_0: \P^y(X_{n_\da+1} \in M_1|X_{n_\da} \in M_0) \ge c_\da 
. 
\end{align}
\item[(c)] 
For $ c_\al, c_\da, n_\da$ as in (b) and for $N := \max \Big\{n_\da,
n_\nu \Big\lceil\frac{\ln (\frac{2- c_\da }{c_\da c_\al})}{-\ln(1-c_\nu c_\nu')} \Big\rceil \Big\}$ 
we have 
\begin{equation}
\label{equ:onsetmonotone}
\forall n \ge N, \forall x,y \in M_0: 
\P^y(X_{n+1} = x) \le \P^y(X_n = x). 
\end{equation}
\end{itemize}
\end{Thm}

\begin{Rem} Onset of monotonicity for transition probabilities of Markov chains. 
\begin{itemize}
\item The condition $|M_0| > 1$ excludes trivial cases. In case of $M_0 = \{x\}$ 
we note that $\P^{x}(X_n = x) = \P^x(X_1 =x)^n$. 
\item 
We note that  
(a) is a result of \cite{CV}, in a formulation for Markov chains in discrete time, 
(b) is easily implied by the assumptions on the chain, and (c) follows from (a) and (b).
\item 
In order to apply (c) and get a concrete value for $N$, we need 
to give concrete values for all the control parameters 
$\nu, n_\nu, c_\nu, c_\nu', c_\al, c_\da, n_\da$. 
We will see how to do that in applying the above criterion 
to the percolation problem we are interested in. 
\end{itemize}
\end{Rem}

We now turn to our application in Bernoulli percolation on 
layered graphs.  
Let $G = (V,E)$ be a (simple) graph and $G^\sq$ the corresponding layered graph.  
We introduce some notation 
for vertices and edges in different layers of $G^\sq$. 
Let $V_n := V \times \{n\}$ denote the set of vertices of the $n$-th layer, and let
$$
E_n^h :=  \{\{(u,n),(v,n)\} \!: uv \in E\} \text{ and }   
E_n^v := \{\{(u,n-1),(u,n)\} \!: u \in V\}
$$ 
denote the sets of all horizontal and vertical edges 
of the $n$-th layer respectively and $E_n := E_n^h \cup E_n^v$. 
Thus $G^\sq = (V^\sq,E^\sq)$, where 
$V^\sq = \bigcup_n V_n$ and $E^\sq = \bigcup_n E_n$. 
We write $E_I := \bigcup_{i \in I} E_i$ for $I \subset \Z$. 
Usually $I$ is a discrete interval written in interval notation, 
where an open interval end is interpreted as only 'half' the 
corresponding layer included, e.g. 
$$
E_{[1,3]} = E_1 \cup E_2 \cup E_3, \quad 
E_{(1,3]} = E_1^h \cup E_2 \cup E_3 \quad \text{ or } 
E_{[1,3)} = E_1 \cup E_2 \cup E_3^v. 
$$
We have suppressed the dependency on $G$ in all the above notations. 
For some fixed vertex $o \in V$ we consider $(o,0) \in V^\sq$
to be the origin of an infection. 

In the following we investigate the percolation process on $G^\sq$ 
via a Markov chain $\cX_n$ that was introduced 
in \cite{KR}. For the sake of completeness we 
properly define $\cX_n$, but we refer to \cite{KR} for further details. 
Informally $\cX_n$ describes, which vertices in layer 
$V_n$ are connected to each other 
and which are connected to $(o,0)$ (and thus infected), 
where in both cases connections constitute of open paths 
that are fully contained in $E_{(-\infty,n]}$, see Figure 
\ref{Fig:Markov} for a visualization.  
The above is formalized in the following definition: 
\begin{figure}[htb!] 
\centering
\vspace*{0.1 cm}
\begin{tikzpicture}[scale = 0.6]
\draw[dotted] (0,-1) grid (4.9,3);
\draw[fill] (0,0) circle (3 pt);
\draw[thick] 
(1,0)--(0,0)--(0,1)--(0,2)--(1,2)--(2,2)--(2,1)
(3,1)--(3,0)
(0,-1)--(1,-1)--(2,-1)--(2,0)
(5,-1)--(4,-1)--(4,0)
(4,1)--(5,1)
(3,2)--(4,2)--(4,3)--(5,3)
(1,3)--(2,3)
;
\draw (6,0) node[anchor=west] {$\cX_0 = \{\{*,0,1\},\{2,4\}, \{3\}\}$};
\draw (6,1) node[anchor=west] {$\cX_1 = \{\{*,0,4\},\{1\},\{2\},\{3\}\}$};
\draw (6,2) node[anchor=west] {$\cX_2 = \{\{*,0,1,2\}, \{3,4\}\}$};
\draw (6,3) node[anchor=west] {$\cX_3 = \{\{*\},\{0,4\},\{1,2\},\{3\}\}$};
\end{tikzpicture}
\caption{The picture shows a realization of Bernoulli percolation on $C_5^\sq$ restricted to $C_5 \times \{-1,0,1,2,3\}$.
The vertices on the right side should be identified with the 
corresponding vertices on the left side,  
the origin $(0,0)$ is marked with a thick dot, open edges are drawn, closed edges are dotted. 
The corresponding values of the chain $\cX_n$ are given
on the right. 
}
\label{Fig:Markov}
\end{figure}
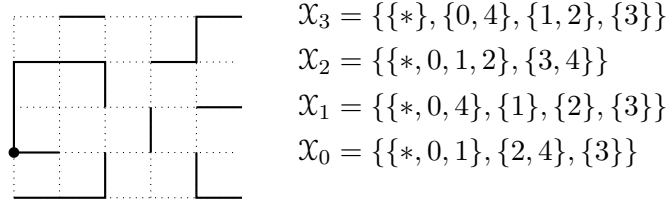
\begin{Def} We consider percolation on $G^\sq$ for some 
graph $G = (V,E)$, 
and fix $o' = (o,0)$ for some $o \in V$. 
\begin{itemize}
\item[(a)]
Let $* \notin V$ be an abstract symbol. 
A partition $x$ of $V \cup \{*\}$ will be called a pattern on $V$. 
Let $\sim_x$ denote the corresponding equivalence relation on $V \cup \{*\}$. 
Let $M$ denote the set of all patterns on $V$ and let
$$
M^* := \{x \in M: \exists u \in V: u \sim_x *\} \quad \text{ and } \quad  
\Mda := M \setminus M^*. 
$$
\item[(b)]  
For $n \ge 0$  and $u',v' \in V^\sq$ let 
$\{u' \lra_{\le n} v'\}$ denote the event that 
there is path from $u'$ to $v'$ consisting of open bonds from $E_{(-\infty,n]}$. 
\item[(c)] 
For $n \ge 0$ and $x \in M$ let $\{\cX_n = x\}$ 
denote the event that 
\begin{align*}
&\forall u,v \in V \!: u \sim_x v  \Lra (u,n) \lra_{\le n} (v,n)
\text{ and } \\
&\forall u \in V: u \sim_x * \Lra  (u,n) \lra_{\le n} o'. 
\end{align*}
\item[(d)] 
For $n \ge 0$ let $\cW_n \!:= |\{v' \! \in V_n \!\!: v' \! \lra \! o'\}|$.
\end{itemize}
\end{Def}

In case of finite $G$, it is intuitively clear that 
$\cX_n$, $n \ge 0$, is a Markov chain, 
since for every $n$ the pattern in the layer $n$ is completely determined 
by the pattern in layer $n-1$ and the bond configuration in $E_n$, 
see Figure \ref{Fig:Markov}.  
More details are given in Section \ref{Sec:Markov}. 
We are now ready to state our main result, 
which is formulated for cycle graphs $G = C_k$ with $k$ vertices, 
where $V = \{0,...,k-1\}$ and $E = \{\{i,i+1\}: i \in V\}$, 
identifying $0$ with $k$.

\begin{Thm} \label{Thm:pattern} 
We consider Bernoulli percolation on $C_k^\sq$ with parameter $p$ for $k \ge 3$. 
Let $o \in C_k$ and $N(k) := 500 k^6 1.95^k$. 
For all $n \ge N(k)$ and $p \in (0,1)$  
\begin{align}
\label{equ:monopattern} 
\forall x \in M^*:~ &\P_p(\cX_n = x) \ge \P_p(\cX_{n+1} = x), \quad \text{ and thus }\\
\label{equ:monoconnection} 
\forall v \in C_k: ~ &\P_p((o,0) \lra (v,n)) \ge 
\P_p((o,0) \lra (v,n+1)) \quad \text{and } \\
\label{equ:monoexpectation}
&\E_p(\cW_n) \ge \E_p(\cW_{n+1}).
\end{align}
\end{Thm}

\begin{Rem}
The onset of monotonicity for Bernoulli percolation on $C_k^\sq$. 
\begin{itemize}
\item 
While \eqref{equ:monoconnection} and thus \eqref{equ:monoexpectation} 
are conjectured to hold for all $n \ge 0$ (even if it is unclear how this could be proved), \eqref{equ:monopattern} clearly does not, as can be 
seen by considering a pattern, for which $o$ is not infected. 
Whereas \cite{KR} shows the existence of a threshold $N(k)$ as above, 
the point of our result is to show by which 
means an explicit value of $N(k)$ can be found. 
\item 
Our value for $N(k)$ grows exponential in $k$. 
While this may seem to be excessively large in terms of the size $k$ of the underlying graph $C_k$, 
it may be more fair to compare this to the size of the state space of the Markov chain $\cX_n$, which is also exponential in $k$. 
\item 
We stress that we do not claim that the given value of $N(k)$
is as small as possible. 
Indeed the CAS-assisted computations 
outlined in \cite{KR} show that our result still is true 
choosing $N(2) := 2, N(3) := 2, N(4) := 4$, 
and we believe that determining the minimal value for $N(k)$ 
for larger $k$ and the order of magnitude for $k \to \infty$ 
is a difficult problem. 
We also note that we have put considerable effort into 
optimizing our method, so we believe that significant 
improvements of the given value of $N(k)$ 
wold require some new ideas of a completely different approach. 

%
%
%
\item 
One may also want to find an explicit value of $N(k,p)$ for fixed $k \ge 3$ and $p \in (0,1)$ so that \eqref{equ:monopattern} holds for all 
$n \ge N(k,p)$. Such a value of $N(k,p)$ may also be extracted from our proof. We note that for a large range of values of $p$, namely $p \in [0.3,0.5]$, the value for $N(k,p)$ will not be significantly smaller than the given value of $N(k)$, whereas for values of $p$ close to $0$ or close to $1$ we get something significantly smaller.  
\item 
For many finite graphs $G$ our proof can be adapted to 
give an explicit value of $N(G)$, such that \eqref{equ:monopattern}
holds for Bernoulli percolation on $G^\sq$ for arbitrary $p \in (0,1)$
and $n \ge N(G)$. 
The condition on $G$ for this to work is that 
the set of attainable states $\{x \in M^*: \exists n \ge 0: \P_p(\cX_n = x) > 0\}$ is communicating. 
For further comments on that see Section 4. 
\end{itemize}
\end{Rem}

\bigskip 

As in \cite{KR} we would like to complement the above result by showing
that \eqref{equ:monoexpectation} also holds uniformly in $n$ for small values of $p$ by a more direct recursive approach. 
Again we are interested in an explicit bound on $p$ for this to be true
in case of $G = C_k$ (or indeed for $G = \Z$). 
\begin{Thm} \label{Thm:numberZ2}
Let $G = C_k$ or $G = \Z$. 
For percolation on $G^\sq$ we have 
\begin{align} \label{equ:numberZ2}
&\forall n \ge 0, p \in [0,0.35]: \E_p(\cW_n) \ge \E_p(\cW_{n+1}).
\end{align}
\end{Thm}

\begin{Rem} \label{Rem:number}
Monotonicity of the expected number of infected vertices. 
\begin{itemize}
\item 
For the graphs considered here this improves the more general result 
of \cite{KR}. The proof of the above theorem gives some ideas 
on how to improve this general result for other specific graphs. 
\item 
Further monotonicity results on the expected number of vertices 
are given in \cite{Ri1}. For $G = C_k$ or $G = \Z$ and 
for all for all $n \ge 0, p \in [0,1]$ we have 
$\E_p(\cW_n) \ge \E_p(\cW_0)$ and we have 
$\E_p(\cW_n) = \infty$ iff  $\E_p(\cW_0) = \infty$. 
\end{itemize}
\end{Rem}

\section{Monotonicity of transition probabilities}
\label{Sec:quasi}

Here we prove Theorem \ref{Thm:quasi}. 
For (a) we note that $\P^y(X_n \in M_0) > 0$ for all $y \in M_0$ and $n \ge 0$ 
(since $M_0$ is communicating with $|M_0| > 1$ and $M_1$ is absorbing), 
 and $\P^y(X_n \in M_1) \to 1$ for $n \to \infty$ and all $y \in M_0$
 (since $M_0$ is transient and finite, and $M_1$ is absorbing). 
 The claim follows by the proof given in Section 5.1 of \cite{CV}. 
While the cited proof uses the setting of Markov chains in continuous time, the proof carries over to our setting of Markov chains in discrete time without change.  

(b) is an easy consequence of the assumptions on the chain: 
The existence of $\al$ is guaranteed by (a) and since $\al$ is 
a distribution we have $\al(y) > 0$ for some $y \in M_0$.
By quasi-stationarity of $\al$ for any $x \in M_0$ we have $\al(x) \ge \al(x) \P^y(X_n = x)$, 
and since $M_0$ is communicating this implies $\al(x) > 0$. 
Since $M_0$ is finite we can choose an appropriate $c_\al > 0$. 
For the second assertion we note that 
there is a constant $\la \in (0,1)$ such that 
$\P^\al(X_n \in M_0) = \la^n$ for all $n \ge0$. 
Propopsition 1.2 of \cite{CV} gives $\P^y(X_n \in M_0) \sim c_y \la^n$ 
for $n \to \infty$ for every $y \in M_0$, which implies $\P^y(X_{n+1} \in M_0 | X_n \in M_0) \to \la$. Since $M_0$ is finite, we can choose appropriate $c_\da \in (0,1)$ and $n_\da \ge 0$.

(c) follows from (a) and (b). 
Consider suitable 
$\nu, n_\nu,c_\nu,c_\nu',c_\al, c_\da, n_\da$ 
and let $c := c_\nu c_\nu'  \in (0,1]$ and $N$ as given. 
By \eqref{equ:cda} for all $y \in M_0$ and $n \ge n_\da$ we have 
\begin{align*}
&\P^y(X_{n+1} \in M_0) = \sum_{x \in M_0} 
\P^y(X_{n-n_\da} = x)\P^x(X_{n_\da+1} \in M_0)\\
&\le  (1-c_\da) \sum_{x \in M_0} 
\P^y(X_{n-n_\da} = x)\P^x(X_{n_\da} \in M_0)
=  (1-c_\da) \P^y(X_{n} \in M_0)
\end{align*}
By (a) we have a quasi-stationary distribution $\al$ such that for all $x,y \in M_0$ and $n \ge 0$ 
\begin{equation*} 
|\P^y(X_n = x|X_n \in M_0) - \al(x)| \le (1-c)^{\lfloor \frac n {n_\nu} \rfloor}
=: \ga_n, 
\end{equation*}
where we have discarded a factor of $2$ by well known properties of the 
total variation distance between probability measures. 
By choice of $N$ for all $n \ge N$ we have $\ga_n \le \frac{c_\al c_\da}{2 - c_\da}\le c_\al$
and for all $n \ge N$ and  $x,y \in M_0$ we have 
\begin{align*}
&\frac{\P^y(X_{n+1} = x)}{\P^y(X_n =x)}
= \frac{\P^y(X_{n+1} = x|X_{n+1} \in M_0)\P^y(X_{n+1} \in M_0)}{\P^y(X_n = x|X_n \in M_0)\P^y(X_n\in  M_0)}\\
&\le  \frac{\al(x) +   \ga_{n+1}}{\al(x)  -  \ga_n} 
(1- c_\da)
\le \frac{c_\al + \ga_n}{c_\al - \ga_n} (1- c_\da)  \le 1, 
\end{align*}
where in the second step we have used the above estimates
noting that $\ga_n \le c_\al \le \al(x)$ by choice of $n$, 
in the third step we have used $\ga_{n+1} \le \ga_n$
and $\al(x) \ge c_\al$, and the last step is equivalent to 
$\ga_n \le \frac{c_\al c_\da}{2 - c_\da}$. \qed

\bigskip 

\section{The Markov chain of patterns} \label{Sec:Markov}

In this section we collect some properties of the Markov chain of patterns. The following definition 
and the following auxiliary results are taken from \cite{KR}. 
We first describe a generalized version of $\cX_n$, denoted by
$\cX_n^{y,k}$, which describes the pattern at layer $n$ that is produced 
by percolation on $E_{[k+1,n]}$ and a prescribed pattern $y$ at layer $k$. 
More formally: 

\begin{Def} We consider percolation on $G^\sq$. 
For $0 \le k \le n$, $x,y \in M$, and $u',v' \in V \times \{k,...,n\}$
let 
\begin{align*}
&E_{k,y} := \{(v,k)(u,k): v,u \in V \!\text{ s.t. }\! v \sim_y u\}, \; 
O_{k,y} := \{(v,k): v \in V \!\text{ s.t. }\! v \sim_y *\},\\  
&\{u' \lra_{y,k,n} v'\} := \{ \exists m \ge 0, v_0',...,v_m' \in V^\sq: 
 v_0' = u', v_m' = v' \text{ and } \\
 &\hspace{3 cm} \quad \forall i:  (v_i'v_{i+1}' \in E_{[k+1,n]}, \cZ_{v_i'v_{i+1}'} = 1) \text{ or } v_i'v_{i+1}' \in E_{k,y}\} \; \text{ and } \\
&\{\cX^{y,k}_{n} = x\} := \{
\forall u,v \in V: u \sim_x v \Lra (u,n) \lra_{y,k,n} (v,n) \text{ and } \\ 
&\hspace{ 3 cm} \forall u \in V: u \sim_x * \Lra \exists w' \in O_{k,y}: 
(u,n)  \lra_{y,k,n} w'\}. 
\end{align*}
Finally, let $\cX^y_n := \cX^{y,0}_n$. 
\end{Def} 
We note that $\cX_n^{y,n} = y$ and in particular 
$\cX_0^y  = y$.  
\begin{Lem} 
We consider percolation on $G^\sq$ with parameter $p \in (0,1)$. 
For all $n \ge 0$ and $y,x_0,..,x_n \in M$ 
we have 
\begin{align*} 
&\P_p(\forall 0 \le i \le n: \cX_i = x_i) = 
\P_p(\cX_0 = x_0) \prod_{1 \le i \le n} \P_p(\cX_1^{x_{i-1}} = x_i) \text{ and }\\
&\P_p(\forall 0 \le i \le n: \cX_i^y = x_i) = 
\de_y(x_0) \prod_{1 \le i \le n} \P_p(\cX_1^{x_{i-1}} = x_i). 
\end{align*} 
In particular $(\cX_n)_{n \ge 0}$ and 
$(\cX_n^y)_{n \ge 0}$ are Markov chains with the same transition probabilities $\pi_p(x,x') := \P_p(\cX_1^x = x')$. 
\end{Lem}
The proof of the lemma is straightforward and can be found in \cite{KR} (Lemma~1). 
We next collect simple properties of the Markov chain of patterns. 
Let the pattern without infection and connections and the pattern, where
everything is connected and infected, be denoted by 
\begin{equation*}
\xda := \{\{*\}\} \cup \{\{v\}: v \in V\} \quad \text{ and } \quad 
\xst := \{\{*\} \cup V\}.
\end{equation*}

\begin{Lem} \label{lem:transitions}
For the above Markov chain of patterns we have: 
\begin{align} 
\label{equ:positivity}
&\forall p,p' \in (0,1), x,y \in M: \pi_p(x,y) > 0 \Lra \pi_{p'}(x,y) > 0. \\
\label{equ:absorbing}
&
\forall p \in (0,1), x \in \Mda, y \in M^*: \pi_p(x,y) = 0. \\
\label{equ:aperiodic}
&\forall p \in (0,1),x \in M: \pi_p(x,x) > 0. \\
\label{equ:nobonds}
&\forall p \in (0,1),x \in M: \pi_p(x,\xda) > 0. 
\end{align}
\end{Lem}
These properties are very easy to see, cf. Lemma 2 of \cite{KR}. 

\bigskip

When investigating the Markov chain $\cX_n$, we restrict our
attention to patterns that actually are attainable in that they 
occur with positive probability:

\begin{Def} For a given graph $G$ and a fixed origin of infection $o \in V$ 
let 
$$
\tM^* := \{x \in M^*: 
\exists n \ge 0: \P_p(\cX_n = x) > 0\}
$$ 
denote the set of all infected attainable states 
of the chain $\cX_n$. 
\end{Def}
By \eqref{equ:positivity}
the positivity of $\P_p(\cX_n = x)$ does not depend on the 
value of $p \in (0,1)$. In case of $G = C_k$ the attainable states can easily be characterized.

\begin{Prop} \label{Prop:Ckcomm}
Consider $G = C_k$ for $k \ge 3$. 
\begin{itemize}
\item[(a)] We have $x \in \tM^*$ iff
$x \in M^*$ and $x$ is noncrossing in that 
$$
\forall v_1 < v_2 < v_3 < v_4 \in \{0,...,k-1\}: \quad 
v_1 \sim_x v_3, v_2 \sim_x v_4 \Ra v_1 \sim_x v_2. 
$$
\item[(b)]
Let $n_k := \lfloor \frac{k+2} 2 \rfloor$. 
For all  $x,y \in \tM^*$ and $p \in (0,1)$  we have   
\begin{equation} \label{equ:estimatecomm}
\P_p(\cX_{n_k}^{y} = x)  \ge  p^{\frac{k^2+2} 2} (1-p)^{\frac{k^2+k} 2}. 
\end{equation}
\end{itemize}
\end{Prop}

\Pf We will first prove two partial results (a') and (b'): 

(a') Every $x \in \tM^*$ is noncrossing: 
Suppose that $\P_p(\cX_n = x) > 0$ and consider a bond configuration 
contributing to $\{\cX_n = x\}$. 
If $v_1 \sim_x v_3$ and $v_2 \sim_x v_4$, then there is an open path in $E_{(-\infty,n]}$ joining $(v_1,n)$ and $(v_3,n)$ 
and another one joining  $(v_2,n)$ and $(v_4,n)$. 
Due to the geometry of $C_k^\sq$ these paths have to meet at some vertex
(formally, this can be seen by a discrete version of the Jordan curve theorem). Thus all four vertices (and in particular $(v_1,n)$ and $(v_2,n)$) are part of the same open cluster, 
which implies $v_1 \sim_x v_2$. 

(b') \eqref{equ:estimatecomm} holds for all $p \in (0,1)$ 
and all $x,y \in M^*$ such that $x$ is noncrossing: 
We fix such $x,y,p$. For $x = x_*$ we can go straight to $V_{n_k}$ 
from an infected vertex of $I$ and then infected all of $V_{n_k}$, 
which gives $\P_p(\cX_{n_k}^{y} = x_*) \ge p^{n_k +k-1} \ge p^{\frac 3 2  k}$ and we are done. Thus we may assume that $x \neq x_*$ and 
w.l.o.g by symmetry of $C_k$ we may assume that $0 \sim_x *$ and 
$k-1 \not \sim_x *$. We first consider the case $0 \sim_y *$. 
We will construct a bond configuration that turns 
$y$ into $x$ in $m$ steps. 
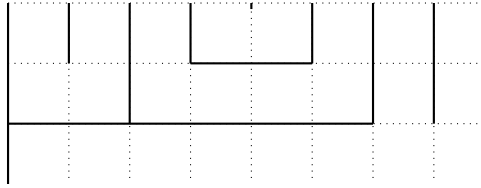
\begin{figure}[htb!] 
\centering
\vspace*{0.5 cm}
\begin{tikzpicture}[scale = 0.8]
\draw[dotted] (0,1.1) grid (7.9,4);
\draw[thick](0,1)-- (0,4);
\draw[thick](1,3)-- (1,4);
\draw[thick](2,2)-- (2,4);
\draw[thick](3,3)-- (3,4);
\draw[thick](5,3)-- (5,4);
\draw[thick](6,2)-- (6,4);
\draw[thick](7,2)-- (7,4);
\draw[thick](4,3.9)-- (4,4);
\draw[thick](0,2)-- (6,2);
\draw[thick](3,3)-- (5,3);
\end{tikzpicture}
\caption{For $C_k$ with $k = 8$, this shows a bond configuration in $E_{[1,3]}$ 
contributing to $\{\cX_3^{y} = x\}$, where $x = \{\{*,0,2,6\}, \{1\}, \{3,5\}, \{4\}, \{7\}\}$ and 
$y \in M^*$ such that $0 \sim_y *$.  
The top layer is $V_3$, drawn bonds are open, dotted bonds are closed, and the right side should be thought of identified with the left side. 
We note that in the given construction  
$C_1^1 = \{0,2,6\}$, $C_1^2 = \{7\}$, $B_1^1 = \{0,...,6\}$, $B_1^2 = \{7\}$, 
$A_1 = \{1,3,4,5\}$, $x_1 = \{\{*,0,1,2,3,4,5,6\},\{7\}\}$, 
$C_2^1 = \{1\}$, $C_2^2 = \{3,5\}$, $B_2^1 = \{1\}$, $B_2^2 = \{3,4,5\}$, 
$A_2 = \{4\}$, $x_2 = \{\{*,0,2,6\}, \{1\}, \{3,4,5\},\{7\}\}$, 
$C_3^1 = B_3^1 = \{4\}$, $A_3 = \emptyset$, $x_3 = x$. 
Exactly $24 - 4$ of the 48 bonds are open.  
The partition $x$ consists of $l= 5$ components and has a nesting number of $m = 3$. 
}
\label{Fig:Ckcomm}
\end{figure}
The following formal description of this construction is illustrated  in Figure \ref{Fig:Ckcomm}. 
For $A \subset \{0,...,k-1\}$ we call the smallest discrete interval containing $A$ the span of $A$. 
We proceed by constructing a suitable bond configuration layer by layer. 
For the first step we choose components $C_1^1,...,C_1^{k_1}$ 
of $x$ such that their spans $B_1^1,...,B_1^{k_1}$ 
form a disjoint decomposition of $\{0,..,k-1\}$
and $C_1^1$ is the infected component of $x$. 
(This is possible since $x$ is noncrossing.) 
Let $A_1 = \{0,...,k-1\} \stm (C_1^1 \cup ... \cup C_1^{k_1})$
and let $x_1$ be the pattern with components $B_1^1,...,B_1^{k_1}$, 
where $B_1^1$ is the infected component of $x_1$. 
$\{X_1^{y,0} = x_1\}$ can be achieved by using one vertical open bond to infect the infected component and horizontal open bonds to join the vertices of each $B_1^i$. 
If $A_1 = \emptyset$ then $x_1 = x$ and we are done. 
Otherwise, in the second step, we choose components 
$C_2^1,...,C_2^{k_2}$ 
of $x$ such that their spans $B_2^1,...,B_2^{k_2}$ 
form a disjoint decomposition of $A_1$. 
(This is possible since $x$ is noncrossing.)   
Let $A_2 = A_1 \stm (C_2^1 \cup ... \cup C_2^{k_2})$
and let $x_2$ be the pattern with components 
$C_1^1,...,C_1^{k_1}, B_2^1,...,B_2^{k_2}$, 
where $C_1^1$ is the infected component of $x_2$. 
$\{X_2^{x_1,1} = x_2\}$ can be achieved by using vertical open bonds to join 
the vertices of $C_1^i$ to the next layer and 
horizontal open bonds to join the vertices of each $B_2^i$. 
If $A_2 = \emptyset$ then $x_2 = x$ and we are done. 
Otherwise we repeat step 2 and carry on in this way. 
If we are done after $m$ steps, we call $m$ the nesting number  of $x$. 
We note that $m \le \lfloor \frac k 2 \rfloor$. 
(The worst case is that 
all but one component of $x$ are nested inside one another, which shows that 
$k \ge 1 + 2(m-1) + 1$.) 
We have thus constructed a bond configuration contributing to 
$\{\cX_m^y = x\}$. It turns out that the number of open and closed 
bonds of this configurations is determined by $k,m$ and the number 
$l$ of components of $x$: 
Identifying the vertices of $V_0$ with the corresponding vertices 
of $V_{m}$ and closing the single open bond in $E_1^v$, 
the remaining open bonds of $E_{[1,m]}$ form a spanning forest 
(on the resulting graph $C_k \times C_m$) consisting of 
$l$ trees. 
By counting the degrees of all vertices of this forest it is 
easy to see that this forest consists of $km - l$ open bonds out of a total of $2km$ bonds. 
%
We thus obtain $\P_p(\cX_m^y= x) \ge p^{km - l+1}(1-p)^{km + l -1}$. 
In case of a general $y \in M^*$ we let $M_0^* := \{y' \in M^*: 0 \sim_{y'} *\}$ and observe that 
$\P_p(\cX_n^y \in M_0^*) \ge p^{n +\lfloor \frac k 2 \rfloor}$ for every $n \ge 1$, since from an infected vertex of $y$ we can go up $n$ steps and then either right or left $\le \lfloor \frac k 2 \rfloor$ steps to get 
to $(0,n)$. Combining the above estimates we obtain 
\begin{align*}
&\P_p(\cX_{n_k}^y \! = x) \ge \P_p(\cX_{n_k}^y \! = x, \cX_{n_k-m} \! \in M_0^*)
= \sum_{y' \in M_0^*} \P_p(\cX_{n_k-m}^y \! =y')\P_p(\cX_m^{y'} \! = x)\\  
&\ge \P_p(\cX_{n_k - m}^y \in M_0^*) \inf_{y' \in M_0^*}\P_p(\cX_m^{y'} = x) 
\ge p^{2n_k + (k-1)m - l}(1-p)^{km + l -1},  
\end{align*}
where we have used the above estimates and $n_k = \lfloor \frac k 2 \rfloor + 1 \ge m+1$. It is easy to see that  $m+1 \le l \le k -(m - 1)$ and thus 
\begin{align*}
&2n_k + (k-1)m - l \le 2 \lfloor \frac k 2 \rfloor + (k-2)m +1  \le k \lfloor \frac k 2 \rfloor + 1 \le \frac {k^2+2} 2 \text{ and }\\
&km + l-1 \le (k-1)m + k \le (k-1)\lfloor \frac k 2\rfloor  + k  \le \frac{k(k+1)} 
2. 
\end{align*}
Using these estimates above we obtain \eqref{equ:estimatecomm}. 

From the above we easily obtain (a) and (b). For (a) one direction is given by (a'), for the other direction let $x \in M^*$ be noncrossing 
and choose some $y \in M^*$ such that $\P_p(\cX_0 = y) > 0$. 
Then (b') implies that $\P_p(\cX_{n_k} = x) \ge \P_p(\cX_0 = y) \P_p(\cX_{n_k}^y = x) > 0$ and thus $x \in \tM^*$. 
(b) directly follows from (a) and (b'). \qed

\begin{Rem} Graphs that can be treated similarly to $C_k$. \nobreak
\begin{itemize}
\item 
We note that the arguments of the following section can be generalized 
to deal with every finite graph $G$ such that for some constant $n_G$ we have 
$$
\forall x,y \in \tM^*, p \in (0,1): \P_p(\cX_{n_G}^{y} = x) > 0. 
$$
In view of \eqref{equ:aperiodic} this is equivalent to the 
property that the attainable infected states form a communicating class w.r.t. the Markov chain $\cX_n$.
\item 
The preceding proposition shows that the above property is satisfied for $C_k$. While it is also easily seen to be satisfied for many 
other graphs, such as complete graphs $K_k$, where $n_{K_k} = 1$, 
it is not satisfied for all graphs: Consider the line graph 
$G = L_3$ with $V = \{0,1,2\}$ and $E = \{\{0,1\}, \{1,2\}\}$ 
with $o := 1$ as the origin of infection.  
Here $x := \{\{1,*\}, \{0,2\}\}$ is attainable, 
since $\P_p(\cX_0 = x) >0$, but $\P_p(\cX_n^{x_*} = x) = 0$ 
for all $n \ge 0$, since a path from $V_0$ to $(1,n)$ makes it impossible for $(0,n), (2,n)$ to be uninfected but connected to each other. 
 \end{itemize}
\end{Rem}

\section{Explicit bounds on the onset of monotonicity} \label{Sec:bounds}

\subsection{The onset of monotonicity}

In this section will prove our main result, Theorem \ref{Thm:pattern}. 
In fact we will prove something slightly stronger: 
\begin{Thm} \label{Thm:patterngeneral}
We consider Bernoulli percolation on $C_k^\sq$ with parameter $p$. 
Let $k \ge 3$, $o \in C_k$ and $N(k) := 500 k^6 1.95^k$. 
For all $n \ge N(k)$ and $p \in (0,1)$  we have 
\begin{equation} \label{equ:monopatterngen}
\forall x,y \in \tM^*: \P_p(\cX_n^y = x) \ge \P_p(\cX_{n+1}^y = x).
\end{equation}
\end{Thm}
In order to see that the above in fact implies Theorem \ref{Thm:pattern} 
it suffices to observe that the monotonicity statements considered are 
ordered in the following way: 
\begin{Prop} \label{Prop:monodiff}
We consider percolation on $G^\sq$ for some finite graph $G$. 
For all $n \ge 0$ and $p \in (0,1)$ 
we have the implications (i) $\Ra$ (ii) $\Ra$ (iii)  $\Ra$ (iv)
of the following statements: 
\begin{enumerate}
\item[(i)] For all $x,y \in \tM^*$ we have $\P_p(\cX_n^y = x) \ge \P_p(\cX_{n+1}^y = x)$. 
\item[(ii)] For all $x \in M^*$ we have 
$\P_p(\cX_n = x) \ge \P_p(\cX_{n+1} = x)$. 
\item[(iii)] 
For all $v \in V$ we have $\P_p((o,0) \lra (v,n)) \ge 
\P_p((o,0) \lra (v,n+1))$.
\item[(iv)]
We have $\E_p(\cW_n) \ge \E_p(\cW_{n+1})$.  
\end{enumerate}
\end{Prop}
For a detailed proof see Proposition 2 and Lemma 3 of \cite{KR}, 
but we will briefly sketch the main ideas: 
(ii) is trivial if $x \in M^*$ is not attainable, and for attainable $x$ 
(ii) can be inferred from (i) by conditioning on the value of $\cX_0$. 
Similarly (iii) can be inferred from (ii) by conditioning on the infection pattern in the layer under consideration. 
Finally,  (iv) follows from (iii) by expressing the number of infected vertices as a sum of indicator variables.  
Preparing for the proof of Theorem \ref{Thm:patterngeneral}, 
we also note the following:

\begin{Prop} \label{Prop:monoN}
We consider percolation on $G^\sq$ for some finite graph $G$. 
Let $p \in (0,1)$ and suppose that for some $N \ge 0$  we have 
$$
\forall y,x \in \tM^*: \P_p(\cX_N^y = x) \ge \P_p(\cX_{N+1} ^y= x),.
$$ 
then \eqref{equ:monopatterngen} holds for all $n \ge N$. 
\end{Prop}
This follows from $\cX_n^y$ being a Markov chain. 
For details see Lemma 3 of \cite{KR}.

The strategy for the proof of Theorem \ref{Thm:patterngeneral} will be to 
use Theorem \ref{Thm:quasi}. Since this works only for values  of $p$ 
bounded away from $0$ and $1$, values of $p$ close to $0$ or $1$ 
have to be treated separately. We start with the main case.  
From now on in this section we will always consider Bernoulli percolation on $C_k^\sq$ for $k \ge 3$ with some parameter $p \in (0,1)$.

\subsection{Intermediate values of $p$}

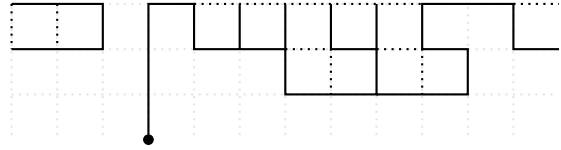
\begin{figure}[htb!] 
\centering
\vspace*{0.5 cm}
\begin{tikzpicture}[scale = 0.6]
\draw[fill] (3,0) circle (3 pt);
\draw[dotted, thick, lgray] (0,0.1) grid (11.9,3);
\draw[thick](3,0)-- (3,3)--(4,3)--(4,2)--(6,2)--(6,1)--(8,1)--(8,3)
(8,1)--(10,1)--(10,2)--(9,2)--(9,3)--(11,3)--(11,2)--(12,2)
(8,2)--(7,2)--(7,3)
(5,2)--(5,3)
(6,2)--(6,3)
(0,2)--(2,2)--(2,3)--(0,3);
\draw[thick, dotted] 
(4,3)--(9,3) (11,3)--(12,3) (6,2)--(7,2)--(7,1) (8,2)--(9,2)--(9,1)
(0,3)--(0,2) (1,3)--(1,2);
\end{tikzpicture}
\caption{The picture shows a realization of Bernoulli percolation on $C_{12}^\sq$ restricted to $E_{[1,3]}$,  
where the dot marks the vertex $(0,0)$, drawn bonds are open, dotted bonds are closed, and the right side should be identified with the left side.
If $0$ is infected with the respect to the pattern $y$, then the shown 
bond configuration contributes to $B \cap B' \subset \{\cX_3^y= x_*\}$. 
}
\label{Fig:plarge}
\end{figure}

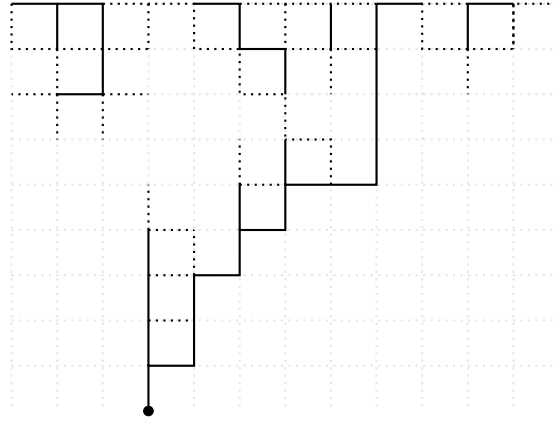
\begin{figure}[htb!] 
\centering
\vspace*{0.5 cm}
\begin{tikzpicture}[scale = 0.6]
\draw[fill] (3,0) circle (3 pt);
\draw[dotted, thick, lgray] (0,0.1) grid (11.9,9);
\draw[thick]
(3,0)--(3,1)--(4,1)--(4,3)--(5,3)--(5,4)--(6,4)--(6,5)--(8,5)--(8,9)
(3,1)--(3,4) (5,4)--(5,5) (6,5)--(6,6)
(8,9)--(9,9) (10,9)--(11,9) (0,9)--(2,9) (4,9)--(5,9)
(1,9)--(1,8) (2,9)--(2,7)--(1,7)  (5,9)--(5,8)--(6,8)--(6,7)
(7,9)--(7,8) (10,9)--(10,8)
;
\draw[thick, dotted] 
(3,2)--(4,2) (3,3)--(4,3) (3,4)--(4,4) (3,4)--(3,5)
(4,3)--(4,4)
(5,5)--(5,6) (5,5)--(6,5)
(7,5)--(7,6)--(6,6)
(0,9)--(0,8)--(1,8)--(3,8) (1,8)--(1,7)--(0,7) (1,7)--(1,6) (2,7)--(2,6)
(2,7)--(3,7)  (3,9)--(3,8) (2,9)--(4,9)--(4,8)--(5,8)--(5,7)--(6,7)--(6,6)
(5,9)--(8,9) (6,9)--(6,8)--(8,8) (7,8)--(7,7)
(10,9)--(9,9)--(9,8)--(11,8)--(11,9)--(12,9) (11,9)--(11,8)
(10,8)--(10,7)
;
\end{tikzpicture}
\caption{The picture shows a realization of Bernoulli percolation on $C_{12}^\sq$ restricted to $E_{[1,9]}$,  
where the dot marks the vertex $(0,0)$, drawn bonds are open, 
dotted bonds are closed, and the right side should be identified with the left side.
If $0$ is infected with the respect to the pattern $y$,   $m' =  5$ and 
$x = \{ \{0\}, \{1,2\}, \{3\}, \{4\}, \{*,5,6\}, \{7,8\}, \{9,10,11\}\}$, 
then the shown bond configuration contributes to $B_1 \cap ... \cap B_4 \subset \{\cX_{m'+4}^y= x\}$. 
%
}
\label{Fig:psmall}
\end{figure}

In this case the main idea is to apply Theorem \ref{Thm:quasi} to $X_n := \cX_n$ and $M_0 := \tM^*$. 
We briefly comment on the intuition on choosing good values for the parameters $\nu,n_\nu,...$. 
The main problem is choosing a good candidate for $\nu$. It is easy to see that for a distribution $\nu$ satisfying \eqref{equ:c1} and \eqref{equ:c2} for given values of $c_\nu,c_\nu'$ any rotation of $\nu$ 
also satisfies these conditions, and it is thus advantageous to choose 
$\nu$ rotation invariant. In that case \eqref{equ:c2} can be seen to hold for $c_\nu' := \frac 1 k$ by symmetry. 
In order to show \eqref{equ:c1} we need to construct suitable bond configurations that allow to go from $y \in \tM^*$ to $x \in \tM^*$ in 
a fixed number of steps and are not too unlikely as compared to $\nu(x)$. 
It thus makes sense to choose $\nu$ as some distribution of infected 
patterns that arises naturally from Bernoulli percolation. 
One of the simplest such candidates is the distribution 
of patterns induced by percolation on $C_k$, where the infected vertex
is chosen uniformly. With this choice we may use the bonds of $E_{n_\nu}^h$ to make all the connections required by the pattern $x$,
and are left to make sure that we have a suitable infection path from layer $0$ to the appropriate vertex of layer $n_\nu$ and 
we have suitable insulation in that the pattern $x$ formed within 
$E_{n_\nu}^h$ is not destroyed 
by additional connections using the bonds below $E_{n_\nu}^h$. 
(For a graphical representation of this construction, which will be used to obtain a good estimate in \eqref{equ:c1}, see Figure \ref{Fig:psmall}.) 
$n_\nu$ should be chosen so that the infection from a given point of 
$V_0$ is sufficiently smeared out over $V_{n_\nu}$, which motivates a choice of the form $n_\nu = c k$. 
The insulation referred to above is too costly in case of larger values of $p$ and it turns out that in this case $\nu = \de_{x_*}$ 
is a better choice, so this will have to be considered additionally. (For a graphical representation 
of the construction used in this case to improve on the estimate in \eqref{equ:c1}, see Figure \ref{Fig:plarge}.) 
One might think that similarly in case of very small values of $p$ 
it may be advantageous to choose $\nu$ as the uniform distribution 
on the $k$ infected totally disconnected patterns, but we found this not to be the case. The choices of $c_\al, c_\da$ are unproblematic. 
We will obtain \eqref{equ:cal} from the estimate of $n_k$-step transition probabilities in the last section, 
and in \eqref{equ:cda} we may simply set $n_\da = 0$ 
since in any step the infection has a decent chance to die out. 
Before carrying out the above program in the proposition below we introduce some constants needed for its formulation: For $p \in (0,1)$ let 
\begin{align}
&c_1(p) := p + (1-p)p^3+ 3(1-p)^2p^5 + 9 (1-p)^3p^7,  \label{equ:c1p} \\
&c_2(p) := 
p+2p^2+2p^3+2p^4+4p^5+8p^6+\frac{64}{1-3p} p^7, 
 \label{equ:c2p} \\
&c_3(p) := \Big(\sqrt{\frac 1 {p^2} + \frac 1 4 e^2 } - \frac 1 2 e \Big)^{-1},
 \label{equ:c3p}   
\end{align}
where we interpret $c_2(p) = \infty$ in case of $p \ge \frac 1 3$. 
We have chosen these constants so that 
$c_i(p) \sim p$ for $p \to 0+$ for all $i$.


\begin{Prop}
\label{Prop:estintermediate}
Let $m := \frac k {2ec_3(p)}$. 
Then \eqref{equ:monopatterngen} is satisfied in the following cases: 
\begin{itemize}
\item[(a)]  $\displaystyle n \!\ge\! N_1(k,p) :=  \!
 \frac{-\ln(p(1 \!-\! p))  (m \!+\! 5) (k^2 \!+\! \frac 3 2 k \!+ \! 2) k^{5/2}  c_2(p)^{m+4}}{
p^{m+4} (\frac{1-p}{1-p(1-p)})^{\frac k 2} (1-p) c_1(1-p)^{k-2}}$, \;  $p \!\in\! (0,\frac 1 3)$,  
\item[(b)]  $\displaystyle n \ge N_2(k,p) :=  
 \frac{- \ln(p(1-p))  (m+5)(k^2+ \frac 3 2 k+2) k^{3/2}}{
p^{m+5} (\frac{1-p}{1-p(1-p)})^{m+1} (1-p) c_1(1-p)^{k-2}}$, \; $p \in (0,\frac 1 2]$,   
\item[(c)] 
$\displaystyle n \ge N_3(k,p) :=   \frac {-  \ln(p(1-p))
\frac 3 2 (k^2+3 k+3)}{p^4  c_1(p)^{k-2}}$, \quad 
$p \in (0,1)$.   
\end{itemize}
\end{Prop}

\Pf
We consider the restriction of the Markov chain $\cX_n$ to the states 
$\tM^* \cup M^\da$. This Markov chain has the properties needed to 
apply Theorem \ref{Thm:quasi}: The state space is a subset of $M$ 
and thus finite, by \eqref{equ:absorbing} $M^\da$ is absorbing, 
by \eqref{equ:nobonds} $\tM^*$ is transient, 
and by Proposition \ref{Prop:Ckcomm} $\tM^*$ is communicating and $|\tM^*| > 1$. 
In all cases we choose $n_\da = 0$ and we note that 
$$
n_\nu \Big\lceil\frac{\ln (\frac{2 - c_\da }{c_\da c_\al})}{-\ln(1-c_\nu c_\nu')} \Big\rceil
\le 
n_\nu \Big( 1 + \frac{\ln 2 - \ln (c_\da c_\al)}{c_\nu c_\nu'}\Big)
\le  
n_\nu \frac{3\ln 2 - \ln (c_\da c_\al)}{c_\nu c_\nu'}, 
$$
where we have used $1 \le 2 \ln 2$ and $c_\nu c_\nu' \le 1$ 
in the last step. 
By Theorem \ref{Thm:quasi} the right hand side 
is a bound for the onset of monotoncity for appropriate choices for $n_\nu, c_\nu,c_\nu', c_\da, c_\al$ depending on $k$ and $p$. 
We note that  \eqref{equ:cal} holds for 
$c_\al :=  p^{\frac{k^2+2} 2} (1-p)^{\frac{k^2+k} 2}$, since for all $x \in \tM^*$ we have 
$$
\al_p(x) = \frac{\P_p^{\al_p}(\cX_{n_k} = x)}{\P_p^{\al_p}(\cX_{n_k} \in \tM^*)}
\ge  \sum_{y \in \tM^*}
\al_p(y) \P_p(\cX_{n_k}^y = x) \ge 
p^{\frac{k^2+2} 2} (1-p)^{\frac{k^2+k} 2},
$$
where in the first step we have used the defining property 
\eqref{equ:quasi1} of $\al_p$, in the second step we have 
used $\P_p^{\al_p}(\cX_{n_k} \in \tM^*) \le 1$ and have split 
cases corresponding to the value of $\cX_0^y$, and in the last 
step we have used $\al_p(\tM^*) = 1$ and \eqref{equ:estimatecomm}.
Finally we observe that \eqref{equ:cda} holds for $n_\da = 0$ and 
 $c_\da := (1-p)^k$ (considering a configuration where all vertical bonds in $E_{1}$ are closed). In the above bound we may thus estimate 
$$
3 \ln 2 - \ln(c_\da c_\al) \le 3 \ln 2 - \frac{k^2+2} 2 \ln p - \frac{k^2+3k} 2 \ln (1-p).  
$$
Using  $\ln 2 \le - \frac 1 2 \ln(p(1-p))$ the above can be seen 
to be $\le - \frac{k^2+3k+3} 2 \ln(p(1-p))$. 
In case of $p \le \frac 1 2$ we can do slightly better and use 
$- \ln(1-p) \le - \ln p$ to see that the above is 
$\le -\frac{2k^2+3k+8}4 \ln(p(1-p))$.

For (c) we let $\nu := \de_{x_*}$ denote 
the distribution with full weight on the configuration, where everything is connected and infected. 
We first show \eqref{equ:c1} for $n_\nu := 3$ and  
$c_\nu(p) := p^4 c_1(p)^{k-2}$.
Let $y \in \tM^*$ and let $i$ be an infected vertex w.r.t. $y$. 
By rotation symmetry we may assume that $i = 0$. 
Let $B$ denote the event that all vertices of $V_3$ 
are connected to each other via open paths in 
$E_{(1,3]}$, 
and let $B'$ denote the event that the straight path from $(0,0)$ 
to $(0,3)$ (consisting of 3 bonds) is open, see Figure \ref{Fig:plarge}.  
We have 
$$
\P_p(\cX_3^y = x_{*}|\cX_3^y \in \tM^*)
\ge \P_p(\cX_3^y = x_{*}) \ge \P_p(B \cap B') \ge p^4c_1(p)^{k-2}. 
$$
In the first step we have used $x_* \in \tM^*$, 
in the second step we have used that on $B\cap B'$ the infection at $(0,0)$ 
is transmitted to all vertices of $V_3$, 
and in the third step we have used the FKG inequality, 
$\P_p(B') = p^3$ and estimated $\P_p(B)$ 
by part (a) of the following lemma. 
Since $\nu(x) = 0$ for all $x \neq x_*$ we have thus established \eqref{equ:c1}. 
We next note that \eqref{equ:c2} holds for $c_\nu' := 1$: 
Indeed, for any $y \in \tM^*$ we have 
$\{\cX_n^y \in \tM^*\} \subset \{\cX_n^{x_*} \in \tM^*\}$. 
The conclusion thus follows from 
\begin{align*}
n_\nu \frac{3\ln 2 - \ln (c_\da c_\al)}{c_\nu c_\nu'} \le 
3 \frac{ - \ln (p(1-p)) \frac {k^2+ 3k+3} 2 }{p^4 c_1(p)^{k-2}} = N_3(k,p). 
\end{align*}

For (a) and (b) a suitable choice for $\nu_p$ is the distribution of patterns 
obtained from percolation on $C_k$ with parameter $p$ 
by infecting one of the $k$ vertices at random. 
More formally, let $z \in \{0,1\}^{E(C_k)}$ and $v \in V(C_k)$. 
For $u \in V(C_k)$ let $K_{z}(u)$ denote the connected component of $u$ w.r.t.\!  $z$, i.e.\! the set of all $u' \in V(C_k)$ connected to $u$ 
by paths consisting of bonds $e$ s.t.\!  $z_e = 1$. 
Let $x_z$ be the collection of all such connected components, 
and $x_{z,v}$ be obtained from $x_z$ by replacing the component 
$K_z(v)$ by $K_z(v) \cup \{*\}$.  
Let $\cZ$ denote independent percolation on $C_k$ and 
$\cV$ be a random vertex of $C_k$ chosen uniformly and independently of the percolation. 
Finally, let $\nu_p$ denote the distribution of $x_{\cZ,\cV}$ w.r.t.\! $\P_p$. 
We first show \eqref{equ:c1} for  suitable choices of $c_\nu$ and  
$n_\nu =:  m' +4$, where $m' := \lceil m \rceil \ge 1$. 
Let $y,x \in \tM^*$. 
W.l.o.g.\! we may assume $x = x_{z,v}$ for some $z,v$ as above, 
since otherwise $\nu_p(x) = 0$ and \eqref{equ:c1} is trivially satisfied. 
Let $i \in C_k$ be a vertex infected w.r.t.\! $y$ and let $d \in \{0,...., \lfloor \frac k2 \rfloor \}$ 
denote the graph distance of $i,v$ in $C_k$. 
W.l.o.g.\! (by symmetry) 
we may assume that $i = 0$ and $v = d$. 
Let $B_1$ denote the event that there is an open path from $(0,0)$
to $(d,m'+1)$ in $E_{[1,m'+1]}$, 
let $B_2$ denote the event that the straight path from $(d,m'+1)$ to $(d,m'+4)$ (consisting of 3 bonds) is open, 
let $B_3 = \{x_{\cZ} = x_z\}$, where by abuse of notation 
$\cZ$ is considered as a percolation process on $E_{m'+4}^h$, 
let $B_4$ denote the event that every vertex in $V_{m'+4} \stm \{(d,m'+4)\}$ 
is not connected to any other vertex of $V_{m'+4} \cup \{(d,m'+3),(d,m'+2)\} \cup V_{m'+1}$ 
by an open path in $E_{[m'+2,m'+4)}$. We have 
$\{\cX_{m'+4}^y = x\} \supset B_1 \cap ... \cap B_4$ by construction, 
see Figure \ref{Fig:psmall}. 
Thus 
\begin{align*}
&\P_p(\cX_{m'+4}^y = x_{z,v}) \ge \P_p(B_1 \cap ... \cap B_4) = \P_p(B_1) \cdot ... \cdot \P_p(B_4) \\
&\ge \binom {m'+d} d p^{d+m'+1} (\frac{1-p}{1-p(1-p)})^{m' \wedge d}
 \cdot p^3 \cdot
\nu_p(x_{z,v}) \cdot (1-p)  c_1(1-p)^{k-2},
\end{align*}
where the second step is due to the independence of $B_1,..,B_4$ 
(noting that the events depend on disjoint edge sets), 
and in the last step we have estimated the probabilities $\P_p(B_i)$: 
For $B_4$ we use part (a') of the following lemma. 
For $B_1$ we use part (c) of the following lemma. 
$\P_p(B_2) = p^3$ is obvious. 
For $B_3$ we use 
$$
\nu_p(x_{z,v}) = \P_p(x_{\cZ,\cV} = x_{z,v}) = \P_p(x_\cZ = x_z, \cV \in K_{z,v}) \le  \P_p(x_\cZ = x_z) = \P_p(B_3).   
$$
We next aim at eliminating the dependency on $d$ in the estimate above.  
We claim 
\begin{equation} \label{equ:estbinomial}
\binom{m'+d} d p^d \ge \frac 1 {2 \sqrt k}.
\end{equation}
Since this is trivial for $d = 0$ we may assume that $1 \le d \le \lfloor \frac k 2 \rfloor$. 
Using the following Lemma \ref{Lem:binomial} 
to estimate the binomial coefficient, we obtain 
$$
\binom{m'+d} d p^d \ge \frac 1 {\sqrt{8d}} \Big(ep\sqrt{(\frac {m'} d)^2 + \frac {m'} d}\Big)^d
\ge  \frac 1 {\sqrt{4k}} \Big(ep\sqrt{(\frac {2m'} k)^2 + \frac {2m'} k}\Big)^d.   
$$
Since 
$$
ep\sqrt{(\frac {2m'} k)^2 + \frac {2m'} k} \ge 1 
\quad \Lra \quad \frac {2m'} k \ge \sqrt{\frac 1 {e^2p^2} + \frac 1 4} - \frac 1 2 = \frac 1 {e c_3(p)}, 
$$
which is satisfied by choice of $m'$, we have thus shown \eqref{equ:estbinomial}. Noting that $\frac{1-p}{1-p(1-p)} < 1$ we obtain 
\begin{align*}
&\P_p(\cX_{m'+4}^y = x_{z,v}) \ge  \frac {p^{m'+4}} {2 \sqrt k} (\frac{1-p}{1-p(1-p)})^{m' \wedge \frac k 2}  (1-p)  c_1(1-p)^{k-2} \nu_p(x_{z,v}). 
\end{align*}
Next we would like to estimate $\P_p(\cX_{n}^y \in \tM^*)$. 
Let $\bcW_{n}^{(i)}$ denote the number of vertices of $V_n$ 
connected to $(i,0)$ via an open path in $E_{[1,n)}$. 
For $p < \frac 1 3$ we have 
\begin{align*}
\P_p(\cX_{n}^y \in \tM^*) &\le \P_p(\exists i: \bcW_{n}^{(i)}\! > 0) \le k \P_p(\bcW_{n}^{(0)}\! > 0)
\le k \E_p(\bcW_{n}^{(0)}) \le kpc_2(p)^{n-1}\!\!,  
\end{align*}
where the first step is by definition - $\cX_n^y$ can only be infected by some open path from $V_0$ to $V_n$, 
the second step is by subadditivity and symmetry, 
and in the fourth step we have used part (b) of the following lemma. 
Combining the two previous estimates we obtain  \eqref{equ:c1} choosing 
$$
c_{\nu}(p) := \frac{p^{m'+4} (\frac{1-p}{1-p(1-p)})^{m' \wedge \frac k 2} (1-p) c_1(1-p)^{k-2}}
{2 \sqrt k (1 \wedge kpc_2(p)^{m'+3})}.
$$
We next show \eqref{equ:c2} for $c_\nu' := \frac 1 k$. 
For $v \in  C_k$ let $B_v$ denote the event that 
$(v,0)$ is connected to $V_n$ via an open path in $E_{[1,n]}$. We have 
\begin{align*}
&\P^{\nu_p}_p(\cX_n \in \! \tM^*)
= \sum_x \nu_p(x) \P_p(\cX_n^x \in \! \tM^*)
= \sum_{v,z} \P_p(\cZ = z, \cV = v) \P_p(\cX^{x_{z,v}}_n \! \in \! \tM^*) \\
&\ge  \sum_{v,z} \P_p(\cZ = z)\P_p(\cV = v) \P_p(B_v)
= \sum_v \frac 1 k \P_p(B_v) \ge \frac 1 k \P_p(\cX_n^y \in \tM^*)  
\end{align*}
for any $y \in \tM^*$. Here 
we have used the definition of $\nu_p$ in the second step and 
$\{\cX_n^y \in \tM^*\} \subset \bigcup_v B_v$ in the last step.
The above gives \eqref{equ:c2} for $c_\nu' := \frac 1 k$. 
Thus
\begin{align*}
n_\nu \frac{3\ln 2 - \ln (c_\da c_\al)}{c_\nu c_\nu'}
\le 
(m' + 4) \frac{- \ln(p(1-p))\frac {2k^2+3k+8} 4 2 k^{3/2} (1 \wedge kpc_2(p)^{m'+3})}{
p^{m' + 4} (\frac{1-p}{1-p(1-p)})^{m' \wedge \frac k 2} (1-p) c_1(1-p)^{k-2}}.  
\end{align*}
The latter can be seen to be $\le N_1(k,p)$ 
estimating $1 \wedge kpc_2(p)^{m'+3} \le kpc_2(p)^{m'+3}$, $m' \wedge \frac k 2 \le \frac k 2$
and $m' \le m+1$ (noting that $\frac{c_2(p)}{p} \ge 1$). 
It can also be seen to be $\le N_2(k,p)$ 
estimating  $1 \wedge kpc_2(p)^{m'+3} \le 1$, $m' \wedge \frac k 2 \le m'$ and  $m' \le m+1$.
This finishes the cases (a) and (b). 
\qed

\begin{Lem} 
Let $c_1(p), c_2(p)$ be defined as in  \eqref{equ:c1p}, \eqref{equ:c2p}.  
\begin{itemize}
\item[(a)]  
Let $A$ denote the event that all vertices of $V_3$ are connected to each other via open paths in $E_{(1,3]}$. 
Then we have $\P_p(A) \ge p c_1(p)^{k-2}$. 
\item[(a')]  
Let $A'$ denote the event that every vertex in $V_3 \stm \{(0,3)\}$ 
is not connected to any other vertex of $V_3 \cup \{(0,2),(0,1)\} \cup V_0$ 
by an open path in $E_{[1,3)}$.
Then we have $\P_p(A') \ge (1-p) c_1(1-p)^{k-2}$.

\item[(b)] Let $\bcW_n$ denote the number of vertices of $V_n$ connected to $o' := (0,0)$ via an open path
in $E_{[1,n)}$. 
Then we have $\E_p(\bcW_n) \le pc_2(p)^{n-1}$. 
\item[(c)] 
For $m \ge 0$ and $0 \le d < k$ let $A_{d,m+1}$ denote the event that there is an open 
path from $(0,0)$ to $(d,m+1)$ in $E_{[1,m+1]}$. Then we have 
$\P_p(A_{d,m+1}) \ge \binom{m+d}d  p^{d+m+1} (\frac{1-p}{1-p(1-p)})^{d\wedge m}$. 
\end{itemize}
\end{Lem}

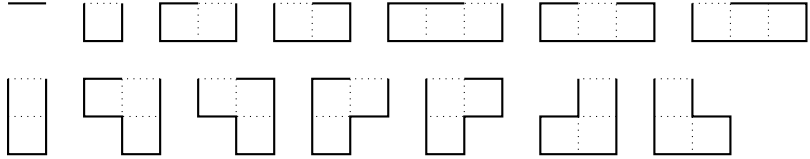
\begin{figure}[htb!] 
\centering
\vspace*{0.5 cm}
\begin{tikzpicture}[scale = 0.5]
\draw[thick] (1,5)--(2,5) 
(3,5)--(3,4)--(4,4)--(4,5)
(6,5)--(5,5)--(5,4)--(7,4)--(7,5)
(8,5)--(8,4)--(10,4)--(10,5)--(9,5)
(13,5)--(11,5)--(11,4)--(14,4)--(14,5)
(16,5)--(15,5)--(15,4)--(18,4)--(18,5)--(17,5)
(19,5)--(19,4)--(22,4)--(22,5)--(20,5)
(1,3)--(1,1)--(2,1)--(2,3)
(4,3)--(3,3)--(3,2)--(4,2)--(4,1)--(5,1)--(5,3)
(6,3)--(6,2)--(7,2)--(7,1)--(8,1)--(8,3)--(7,3)
(10,3)--(9,3)--(9,1)--(10,1)--(10,2)--(11,2)--(11,3)
(12,3)--(12,1)--(13,1)--(13,2)--(14,2)--(14,3)--(13,3)
(16,3)--(16,2)--(15,2)--(15,1)--(17,1)--(17,3)
(18,3)--(18,1)--(20,1)--(20,2)--(19,2)--(19,3)
; 
\draw[dotted] (3,5)--(4,5)
(6,4)--(6,5)--(7,5)
(8,5)--(9,5)--(9,4)
(12,4)--(12,5) (13,4)--(13,5)--(14,5)
(16,4)--(16,5)--(17,5)--(17,4)
(19,5)--(20,5)--(20,4) (21,5)--(21,4)
(1,3)--(2,3) (1,2)--(2,2)
(5,3)--(4,3)--(4,2)--(5,2)
(6,3)--(7,3)--(7,2)--(8,2)
(11,3)--(10,3)--(10,2)--(9,2)
(12,3)--(13,3)--(13,2)--(12,2)
(17,3)--(16,3) (17,2)--(16,2)--(16,1)
(18,3)--(19,3)
(18,2)--(19,2)--(19,1)
;
\end{tikzpicture}
\caption{Connecting horizontally adjacent vertices in $\Z^2$ with shortest open paths of length 1,3,5,7. 
Open bonds are drawn, closed bonds are dotted. 
All paths go at most two layers below (and do not go above) the starting layer.
}
\label{Fig:c0c1}
\end{figure}

\Pf (a) For $G = L_k$ (the line graph on the $k$ vertices $0,...,k-1$) 
let $B_i$ denote the event that the vertices 
$(i-1,3)$ and $(i,3)$ are connected via open paths in 
$E_{(1,3]}$.  We have 
$$
\P_p(A) \ge \P_p(\bigcap_{1 \le i \le k-1} B_i) 
\ge \prod_{1 \le i \le k-1} \P_p(B_i), 
$$
where the first step follows from an obvious coupling, 
the second step follows from the FKG inequality 
(since $B_1,...,B_{k-1}$ are increasing).  
To calculate $\P_p(B_i)$ we consider the shortest 
open path from $(i-1,3)$ to $(i,3)$ in $E_{(1,3]}$ 
and its length $L_i \in \{1,3,5,...\} \cup \{\infty\}$.  
For $3 \le i \le k-3$  we have 
$\P_p(B_i) \ge \P_p(L_i \le 3) = c_1(p)$, where the last step 
follows by counting possible connection paths, see Figure \ref{Fig:c0c1}. 
For $i \in \{1,2,k-2,k-1\}$ we have to be more careful 
in that the paths may not leave $L_k^\sq$. 
Again referring to Figure \ref{Fig:c0c1} we get 
$\P_p(B_i) \ge c_1'(p) = p + (1-p)p^3 + 2 (1-p)p^5$. 
Thus we obtain $\P_p(A) \ge c_1(p)^{k-5} c_1'(p)^4$ (and we note that this 
is still true if $k \in \{3,4\}$, since $c_1(p) \ge c_1'(p)$). The conclusion 
follows from noting that $c_1(p)^{k-5} c_1'(p)^4 \ge c_1(p)^{k-2} p$, for which we refer 
to Lemma \ref{Lem:ac1}.

(a') For $G = L_{k+1}$ let $\tilde A'$ be the event that 
every vertex in $\{1,...,k-1\} \times \{3\}$  is not connected to any other vertex of the inner boundary of  $\{0,...,k\} \times \{0,1,2,3\}$ by an open path in $E_{[1,3)}$. Let $B_i$ as in (a). 
We have 
$$
\P_p(A') = \P_p(\tilde A') \ge \P_{1-p}(\bigcap B_i), 
$$
where the first step follows from an obvious coupling, and the second step follows from the dual coupling, where the edges of $L_k \times \{0,1,2\}$
are identified with the dual edges of the non-boundary edges of 
$L_{k+1} \times \{0,1,2,3\}$, noting that in any configuration contributing to $B_i$ in the dual lattice 
the cluster $(i,3)$ is not connected 
to any other vertex of the inner boundary of $\{0,...,k\} \times \{0,1,2,3\}$ by an open path in $E_{[1,3)}$.
Estimating $\P_{1-p} (\bigcap B_i)$ as in (a)  gives the claim. 

\begin{figure}[htb!] 
\centering
\vspace*{0.5 cm}
\begin{tikzpicture}[scale = 0.6]
\draw[dotted, thick, lgray] (-1,0) grid (18,2);
\draw[thick]
(0,1)--(5,1)--(5,2) 
(6,1)--(8,1)--(8,0)--(9,0)--(9,2)
(10,1)--(11,1)--(11,0)--(12,0)--(12,1)--(13,1)--(13,2)
(14,1)--(15,1)--(15,0)--(17,0)--(17,2)
;
\end{tikzpicture}
\caption{Self-avoiding paths of length 6 in  $E_{(-\infty,1)} \stm E_{o'}$
connecting $o'$ to $V_1$. 
We have listed all such paths starting with a step from $o'$ to the right, 
and we have symmetric paths starting with a step to the left, which shows $n_6 = 8$.
}
\label{Fig:pathcounting1}
\end{figure}
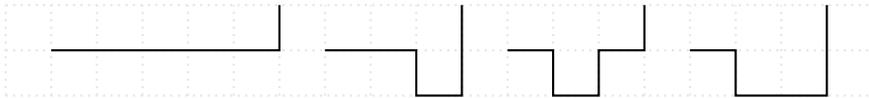
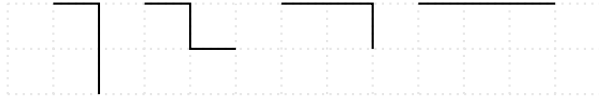
\begin{figure}[htb!] 
\centering
\vspace*{0.5 cm}
\begin{tikzpicture}[scale = 0.6]
\draw[dotted, thick, lgray] (0,0) grid (13,2);
\draw[thick]
(1,2)--(2,2)--(2,0) (3,2)--(4,2)--(4,1)--(5,1) 
(6,2)--(8,2)--(8,1) (9,2)--(12,2) 
;
\end{tikzpicture}
\caption{Self-avoiding paths of length 3 in $E_{(-\infty,0]} \stm E_{o'}$ starting in $o'$ 
with a step to the right. There are symmetric paths starting with a step to the left.  
}
\label{Fig:pathcounting2}
\end{figure}
\begin{figure}[htb!] 
\centering
\vspace*{0.5 cm}
\begin{tikzpicture}[scale = 0.6]
\draw[dotted, thick, lgray] (0,0) grid (12,1);
\draw[thick]
(1,0)--(2,0)--(2,1)--(1,1)
(3,0)--(3,1)--(5,1)
(6,0)--(7,0)--(7,1)--(8,1)
(9,0)--(11,0)--(11,1)
;
\end{tikzpicture}
\caption{Self-avoiding paths of length 3 in $E_{(-\infty,0]}$ connecting $(0,-1)$ to $V_0$. 
For each shown path we also have its mirror image, which shows $k_1 = 8$.
}
\label{Fig:pathcounting3}
\end{figure}
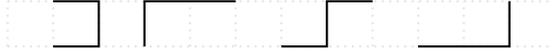

(b)
We proceed inductively. We note that $\E_p(\bcW_1) = p$ since $E_{[1,1)} = E_1^v$. 
For the inductive step we first note that for every $w' \in V_{n+1}$ 
connected to $o'$ by an open path in $E_{[1,n+1)}$, 
the connecting path can be chosen to be vertex-self-avoiding. 
For such a path let $v'$ denote the first vertex in $V_n$. 
Note that the 
preceding vertex has to be the one below $v'$, so the subpath 
from $v'$ to $w'$ is disjoint from $E_{v'}$, 
the edge set consisting of the four edges adjacent to the vertex below $v'$. 
Using a union bound and the BK-inequality (and the $\circ$-notation in its formulation, see e.g. \cite{G}) we thus obtain 
\begin{align*}
&\E_p(\bcW_{n+1}) = \sum_{w' \in V_{n+1}} \P_p(o' \lra w' \text{ w.r.t. } E_{[1,n+1)}) \\
&\le \sum_{w' \in V_{n+1}} \sum_{v' \in V_n} \P_p(\{o' \lra v'  \text{ w.r.t. } E_{[1,n)}\} \circ \{v' \lra w' \text{ w.r.t } E_{(-\infty,n+1)} \stm E_{v'}\})\\
&\le  \sum_{v' \in V_n} \P_p(o' \lra v'  \text{ w.r.t. } E_{[1,n)})
\sum_{w' \in V_{n+1}} \P_p(v' \lra w' \text{ w.r.t } E_{(-\infty,n+1)} \stm E_{v'}).  
\end{align*}
The sum over $w'$ can be written as an expectation and by symmetry it 
equals $\E_p(\bcW)$, 
where $\bcW$ denotes the number of vertices of $V_1$ 
connected to $o'$ via an open path in $E_{(-\infty,1)} \stm E_{o'}$.
Thus we obtain 
$$
\E_p(\bcW_{n+1}) \le \E_p(\bcW_{n}) \E_p(\bcW). 
$$
It remains to check that $\E_p(\bcW) \le c_3(p)$ in case of $p < \frac 1 3$. 
We note that $\E_p(\bcW) \le \sum_{l \ge 1} n_l p^l$, 
where $n_l$ denotes the number of (vertex-)self-avoiding paths 
of length $l$ in $E_{(-\infty,1)} \stm E_{o'}$ connecting $o'$ to $V_1$. 
It is easy to see that 
$n_1 = 1, n_2 = 2, n_3 = 2, n_4 = 2, n_5 = 4, n_6 = 8$, 
for $n_6$ see Figure \ref{Fig:pathcounting1}. 
In order to estimate $n_l$ for $l \ge 7$ we note that a contributing path can be decomposed into subpaths of lengths $3,l-7,3,1$. 
The first subpath is contained in $E_{(-\infty,0]} \stm E_{o'}$
and can be chosen in $8$ ways, see Figure \ref{Fig:pathcounting2}.
For the second subpath we have $\le 3$ choices in every step (since the path cannot backtrace). 
The third subpath has to connect some point $(m,-i)$ to $V_0$. 
Let $k_i$ denote the number of self-avoiding paths of length $3$ starting at $(m,-i)$ and ending in $V_0$. We have $k_0 = 4, k_1 = 8, k_2 = 6, k_3 = 1$ and $k_l = 0$ for $l > 3$, see Figure \ref{Fig:pathcounting3}
for $k_1$.  Thus there are at most 8 choices for this part. 
The last step goes up to layer $V_1$. In conclusion we see 
that $n_l \le 8 \cdot 3^{l-7} \cdot 8 \cdot 1= 64 \cdot 3^{l-7}$, 
which completes the proof of the claim. 

(c) Let $\Ga$ denote the set of all paths from from $(0,0)$ to $(d,m)$ 
of length $d+m$ (say in the graph $\Z^2$). 
We note that every path $\ga \in \Ga$ only consists of steps that go up or to the right. 
Let $A_{d,m}'$ denote the event that some path $\ga \in \Ga$ is open. 
By an obvious coupling we have $\P_p(A_{d,m+1}) \ge p \P_p(A_{d,m}')$. 
We note that there is a partial order on $\Ga$ 
such that $\ga' \preceq \ga$ iff $\ga'$ lies below $\ga$ 
in that for every $0 \le i \le d$ and every vertex $(i,j')$ of $\ga'$ 
there is a vertex $(i,j)$ of $\ga$ such that $j' \le j$. 
Furthermore, for any two paths $\ga,\ga' \in \Ga$ 
there is path $\ga \vee \ga' \in \Ga$ contained in $\ga \cup \ga'$ 
such that $\ga, \ga' \preceq \ga \vee \ga'$. 
Thus the (random) subset of open paths in $\Ga$ contains 
a unique maximal element (w.r.t. $\preceq$). 
For $\ga \in \Ga$ let $A_\ga$ denote the event that $\ga$ is 
this maximal open path. Thus $A_{d,m}'$ is the disjoint union 
of the events $A_\ga, \ga \in \Ga$. 
In order to estimate $\P_p(A_\ga)$, we note that $\ga$ has a vertex set of the form 
$\{(i,j): 0 \le i \le d, j_{i-1} \le j \le j_i\}$ for some sequence $0 = j_{-1} \le j_0 \le ... \le j_d \le m$. 
Let $B_\ga$ denote the event that all edges of $\ga$ are open. 
For  $0 \le i < d$ let $B_{\ga,i}$ denote the event that 
there is no open path starting at $(i,j_i)$, 
going up a number of steps and then going right one step to end in a vertex of the form $(i+1,j)$, 
where $j_i < j \le m$. 
With these notations we have $A_\ga \supset B_\ga \cap \bigcap_i B_{\ga,i}$. (For a typical bond configuration in this intersection see 
the layers $E_{(1,6]}$ of Figure \ref{Fig:psmall}.)  
We note that 
$$
\P_p(B_{\ga,i})
\ge \sum_{k\ge 0} p^k(1-p)^{k+1} = \frac{1-p}{1-p(1-p)}. 
$$
Since the events $B_\ga,B_{\ga,0},..,B_{\ga,d-1}$ depend on disjoint edge sets and thus are independent,
we thus get 
\begin{align*}
&\P_p(A_\ga) \ge \P_p(B_\ga \cap  \bigcap_i B_{\ga,i})  =  p^{d+m} (\frac{1-p}{1-p(1-p)})^d, \text{ implying} \\
&\P_p(A_{d,m}') \ge \P_p(\bigcup_\ga A_\ga) = \sum_\ga \P_p(A_\ga)
\ge \binom{m+d}d  p^{d+m} (\frac{1-p}{1-p(1-p)})^d, 
\end{align*}
where we also have used the disjointness of $A_\ga$, $\ga \in \Ga$, 
and we have counted $|\Ga|$. 
By symmetry the same estimates holds if $m,d$ are interchanged, 
i.e. we also have 
$\P_p(A_{d,m}') \ge \binom{m+d}d  p^{d+m} (\frac{1-p}{1-p(1-p)})^{m}$. 
This gives the desired estimate for $\P_p(A_{d,m+1})$. 
\qed 

\begin{Lem} \label{Lem:binomial}
For all $m, d \ge 1$ we have 
$$
\binom {m+d} d \ge \frac 1 {\sqrt {8}} \sqrt{\frac 1 d + \frac 1 m} \Big( e \sqrt{(\frac{m}{d})^2+\frac{m}{d}}\Big)^d.
$$ 
\end{Lem}

\Pf For $d = 1$ this follows from $\sqrt 8 \le e$. Otherwise we use 
the factorial estimate by Robbins (see \cite{R}):  
$n! = \sqrt{2\pi} \frac{n^{n + \frac 1 2}}{e^n} e^{r_n}$ for all $n \ge 1$, 
where $\frac 1 {12n+1} \le r_n \le \frac 1 {12n}$.  Thus 
\begin{align*}
&\binom {m+d} d =
\sqrt{\frac{m+d}{2\pi m d}} \frac{(m+d)^{m+d}}{m^m d^d} e^{r_{m+d} - r_m - r_d}, \quad \text{ where } \\
&r_m + r_d  - r_{m+d} \le \frac 1 {12 m } + \frac 1 {12 d} - \frac 1 {12(m+d) + 1} =: h(m,d).
\end{align*}
Since  $h$ is decreasing both in $m$ and $d$ and $d \ge 2$, 
we have $h(m,d) \le h(1,2) \le \frac 1 {10}$, and thus 
$\sqrt{\frac{1}{2\pi}}e^{r_{m+d} - r_m - r_d} \ge \sqrt{\frac{1}{2\pi}} e^{-\frac 1 {10}} \ge 
\frac 1 {\sqrt 8}$. Setting $a := \frac m d$ we next note that 
$$
\frac{(m+d)^{m+d}}{m^m d^d} = \Big(\frac{(a+1)^{a+1}}{a^a}\Big)^d = 
\Big((1+\frac 1 a)^{a+\frac 1 2} \sqrt{a^2+a}\Big)^d 
\ge \Big(e \sqrt{a^2 + a}\Big)^d, 
$$
using Lemma \ref{Lem:ae} in the last step.
This gives the desired estimate.  
 \qed

\subsection{Small values of $p$}
 
Our general strategy in this case is similar to the one for general graphs in \cite{KR}: 
If $p$ is small and we want to go from a pattern $y$ to a pattern $x$ in large number of steps, in most layers there will be exactly one open bond. 
In particular, with high probability there will be two consecutive layers with a matching vertical open bond such that discarding 
one of these layers will not destroy pattern $x$. 
This will let us compare  $\P_p(\cX_{n+1}^y = x)$ and 
$\P_p(\cX_{n}^y = x)$. 
Unlike in \cite{KR} we have to do our estimates very careful and 
make use of the simple structure of the graph $C_k$. We will use 
the notation $e_v$ for the (vertical) bond below $v$, i.e. 
$$
e_{(j,l)} := \{(j,l), (j,l-1)\} \text{ for } (j,l) \in C_k \times \Z. 
$$

\begin{Lem} \label{Lem:estsmallp}
For all $n \ge n_k$, $p \in (0,\frac 1 {20}]$ and 
$y,x \in \tM^*$ we have 
 \begin{align} 
\label{equ:estsmallpmatching}
&\P_p(\cX_{n+1}^y = x, A_{n+1}) \le p \frac{n+1} 3 \P_p(\cX_n^y = x),\\
\label{equ:estsmallpnotmatching}
&\P_p(\cX_{n+1}^y = x, A_{n+1}^c) \le 
(kp^4(2 + 10p +2kp))^{\frac {n-1} 3},
\\
\label{equ:estsmallpprob}
&\P_p(\cX_n^y = x) \ge p^{\frac{k^2-k+1}2}(1-p)^{\frac{k^2+k}2} p^{n} ,  
\end{align}
where  $A_{n+1} := \bigcup_{0 < i \le \frac {n+1}3} B_{3i}$, $B_{3i} := \bigcup_{j \in C_k} 
B_{j,3i-2}^{3i} \cap B_{j,3i-3}^{3i}$ and 
$B_{j,l'}^l$ denotes the event that $(j,l')$ is the unique vertex of $V_{l'}$ connected to $V_l$ via an open path in $E_{[l'+1,l)}$. 
\end{Lem}

\Pf 
For \eqref{equ:estsmallpmatching} we use the definition of $A_{n+1}$ and a union bound to obtain
\begin{align*}
&\P_p(\cX_{n+1}^y \!=  x, A_{n+1})
\le \!\! \sum_{i,j,x',y'} \!\!
\P_p(\cX_{3i-3}^y \! = y'\!, B_{j,3i-2}^{3i}, B_{j,3i-3}^{3i},
\cX_{3i}^y \! =  x'\!, \cX_{n+1}^y \! =  x)\\
& \le   \!\! \sum_{i,j,x',y'} \!\!
\P_p(\cX_{3i-3}^y = y')\P_p(B_{j,1}^3, B_{j,0}^3,
\cX_{3}^{y'} = x')\P_p(\cX_{n+1-3i}^{x'} = x)
\end{align*}
where $0 < i \le \frac{n+1} 3$, $j \in C_k$,  $y',x' \in \tM^*$, 
and in the last step we have decomposed the event into three parts depending on disjoint sets of bonds and used independence. 
On $B_{j,1}^3 \cap B_{j,0}^3$ all connections formed in $V_3$ 
via open paths in $E_{(-\infty,3]}$ are in fact already formed via open 
paths in $E_{[2,3]}$ and the infection can only be transmitted from $(j,1)$ and thus from $(j,0)$, which implies  $\cX_3^{y'} = \cX_3^{y',1}$. 
Also on $B_{j,1}^3 \cap B_{j,0}^3$, the bond $e_j := e_{(j,1)}$ is open. 
Thus 
\begin{align*}
&\P_p(B_{j,1}^3, B_{j,0}^3, \cX_{3}^{y'} = x') = 
\P_p(B_{j,1}^3, B_{j,0}^3, \cX_{3}^{y',1} = x') \\
&\ge \P_p(\cZ_{e_{j}} = 1, B_{j,1}^3, \cX_{3}^{y',1} = x') 
=  \P_p(\cZ_{e_{j}} = 1) \P_p(B_{j,0}^2, \cX_{2}^{y'} = x')
\end{align*}
where in the last step we again have decomposed the event into parts depending on disjoint sets of bonds. 
We have $\P_p(\cZ_{e_{j}} = 1) = p$ and since the 
$B_{j,0}^2, j \in C_k$ are disjoint we have 
$\sum_j \P_p(B_{j,0}^2, \cX_{2}^{y'} = x') \le \P_p (\cX_{2}^{y'} = x')$. 
Combining all the above we get  
\begin{align*}
\P_p(\cX_{n+1}^y = x, A_{n+1})
&\le \sum_{i,x',y'}
\P_p(\cX_{3i-3}^y = y')p\P_p(\cX_{2}^{y'} = x')\P_p(\cX_{n+1-3i}^{x'} = x)\\
&= p \sum_i \P_p(\cX_{n}^y = x) \le p \frac{n+1} 3 \P_p(\cX_{n}^{y} = x)
\end{align*}
as claimed. 
For \eqref{equ:estsmallpnotmatching} let 
$B_l'$ be the intersection of $B_l^c$ with the event that 
$V_{l-3}$ and $V_{l}$ are joined by an open path in  $E_{[l-2,l)}$. 
We note that 
$$
\P_p(\cX_{n+1}^y = x, A_{n+1}^c) \le \P_p(\bigcap_{0 < i \le \frac {n+1} 3}B_{3i}')
\le \P_p(B_3')^{\frac{n-1}3}, 
$$
where we have used the definition of $A_{n+1}$, the independence of the $B_{3i}'$ 
(which depend on disjoint bond sets), $\P_p(B_l') = \P_p(B_3')$
and $\lfloor \frac{n+1} 3\rfloor \ge \frac{n-1} 3$. 
With the estimation of $\P_p(B_3')$ 
from the following lemma we are done. 
For \eqref{equ:estsmallpprob} we note that for $n \ge n_k$ we have 
\begin{align*}
&\P_p(\cX_n^y = x) = \sum_{x' \in \tM^*} \P_p(\cX_{n-n_k}^y = x')\P_p 
(\cX_{n_k}^{x'} = x)\\
&\ge \sum_{x' \in \tM^*} \P_p(\cX_{n-n_k}^y = x') p^{\frac{k^2+2} 2}(1-p)^{\frac{k^2+k} 2}  \\
&= \P_p(\cX_{n-n_k}^y \in \tM^*) p^{\frac{k^2+2} 2}(1-p)^{\frac{k^2+k} 2}
\ge p^{n-n_k} p^{\frac{k^2+2} 2}(1-p)^{\frac{k^2+k} 2}, 
\end{align*}
where in the second step we have used \eqref{equ:estimatecomm}, 
and in the last step we have used that an infection 
can be transmitted by a path of open vertical bonds. 
Using $n_k \ge \frac{k+1}2$ gives the desired estimate. \qed 

\begin{Lem} 
Let $B'$ denote the intersection of $B_3^c$ with the event 
that $V_0$ and $V_3$ are connected by an open path in $E_{[1,3)}$.  
For all $p \in (0,\frac 1 {20}]$ we have 
$$
\P_p(B') \le (2+10p+2kp)kp^4. 
$$
\end{Lem}

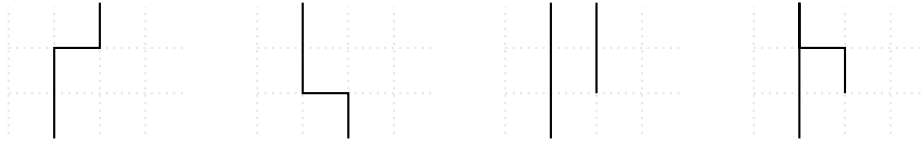
\begin{figure}[htb!] 
\centering
\vspace*{0.5 cm}
\begin{tikzpicture}[scale = 0.6]
\draw[dotted, thick, lgray] (0,0.1) grid (3.9,2.9);
\draw[thick]
(1,0)--(1,2)--(2,2)--(2,3);
\end{tikzpicture} 
\qquad 
\begin{tikzpicture}[scale = 0.6]
\draw[dotted,  thick, lgray] (0,0.1) grid (3.9,2.9);
\draw[thick]
(2,0)--(2,1)--(1,1)--(1,3);
\end{tikzpicture} 
\qquad 
\begin{tikzpicture}[scale = 0.6]
\draw[dotted,  thick, lgray] (0,0.1) grid (3.9,2.9);
\draw[thick]
(1,0)--(1,3) (2,1)--(2,3);
\end{tikzpicture} 
\qquad 
\begin{tikzpicture}[scale = 0.6]
\draw[dotted,  thick, lgray] (0,0.1) grid (3.9,2.9);
\draw[thick]
(1,0)--(1,3) --(1,2)--(2,2)--(2,1);
\end{tikzpicture} 
\caption{Bond configurations in $E_{[1,3)}$ for $k = 4$. 
The right side should be identified with the left side. 
The illustration 
shows typical configurations  connecting $V_0$ to $V_3$ by an open path in $B_3$, $B_\lfloor$, $B_\parallel$ and 
$B_\vdash$ respectively. 
}
\label{Fig:decomposingB}
\end{figure}

\Pf We define a suitable decomposition of $B'$, 
see Figure \ref{Fig:decomposingB} for an illustration of the different cases. 
Let $B_\lfloor$ denote the event that some vertex $v \in V_1$ is connected to $V_3$ by an open path in $E_{[2,3)}$, 
and $v$ is connected to $V_0$ by an open path in $E_{[1,2]}$ disjoint from the first path and disjoint from $e_v$. 
Let $B_\parallel$ denote that event that for some vertices 
$v \neq v' \in V_1$ we have two (bond-)disjoint open paths connecting 
these vertices to $V_3$ in $E_{[2,3)}$ and $e_{v}$ is open.  
Let $B_\vdash$ denote the event that 
for some vertices 
$v \neq v' \in V_1$ we have two (bond-)non-disjoint open paths connecting 
these vertices to $V_3$ in $E_{[2,3)}$ and $e_{v}$ is open.  
We claim
\begin{equation}
\label{equ:inclusion}
B' \subset B_\lfloor \cup B_\parallel \cup B_\vdash.
\end{equation}
On $B'$ we have an open self avoiding path in $E_{[1,3)}$ from $V_0$ to $V_3$. Let $v$ denote the last vertex of the path in $V_1$. 
If the path does not start with $e_v$, we are in $B_\lfloor$. 
On the other hand, if the path starts with $e_v$, by definition of $B_3$
either for some $v' \in V_1$ with $v' \neq v$ we have an open 
path from $v'$ to $V_3$ in $E_{[2,3)}$, 
or for some $v'' \in V_1$ with $v'' \neq v$ we have an open path 
from $v''$ to $V_3$ in $E_{(1,3)}$ and $e_{v''}$ is open. 
In the first case we are in $B_\parallel \cup B_\vdash$. 
In the second case let $v'$ denote the last vertex in $V_1$ 
of this additional path. 
If $v' \neq v$ we are in fact again in the first case, and 
if $v' = v$ we are in $B_\lfloor$. This finishes the proof of \eqref{equ:inclusion}.

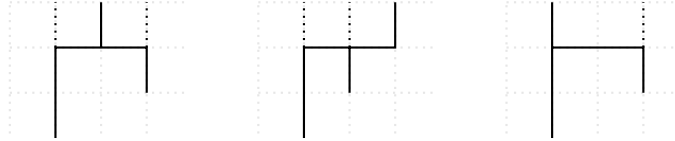
\begin{figure}[htb!] 
\centering
\vspace*{0.5 cm}
\begin{tikzpicture}[scale = 0.6]
\draw[dotted, thick, lgray] (0,0.1) grid (3.9,3);
\draw[thick]
(1,0)--(1,2)--(3,2)--(3,1) (2,2)--(2,3) ;
\draw[thick, dotted]
(1,2)--(1,3) (3,2)--(3,3);
\end{tikzpicture} 
\qquad 
\begin{tikzpicture}[scale = 0.6]
\draw[dotted, thick, lgray] (0,0.1) grid (3.9,3);
\draw[thick]
(1,0)--(1,2)--(2,2)--(2,1) (2,2)--(3,2)--(3,3) ;
\draw[thick, dotted]
(1,2)--(1,3) (2,2)--(2,3);
\end{tikzpicture} 
\qquad 
\begin{tikzpicture}[scale = 0.6]
\draw[dotted, thick, lgray] (0,0.1) grid (3.9,3);
\draw[thick]
(1,0)--(1,3) (1,2)--(3,2)--(3,1) ;
\draw[thick, dotted]
(3,2)--(3,3);
\end{tikzpicture} 

\caption{Bond configurations in $E_{[1,3]}$ for $k = 4$. 
Open edges are drawn and closed edges are dotted, and the right side should be identified with the left side. 
The illustration shows bond configurations corresponding to the two cases within 
$B_\vdash \cap B_\parallel^c$ 
}
\label{Fig:tree}
\end{figure}

For estimating $\P_p(B_\lfloor)$ we first note that 
$B_\lfloor =  \bigcup_{v \in V_1} (B_\lfloor^{v,\uparrow} \circ B_\lfloor^{v,\rightarrow})$, where 
$B_\lfloor^{v,\uparrow}$ denotes the event that $v$ is connected to $V_3$  via a (self-avoiding) open path in $E_{[2,3)}$, 
and $B_\lfloor^{v,\rightarrow}$ denotes the event that $v$ is connected to $V_0$ via a (self-avoiding) open path in $E_{[1,2]}$, which starts with a horizontal step. 
The path from the definition of $B_\lfloor^{v,\uparrow}$ starts with 
one step going up. If the next bond up is closed, then it continues some 
number of steps to the right (or left) and then goes up. Thus 
$$
\P_p(B_\lfloor^{v,\uparrow}) \le p^2 + p(1-p)2(p+p^2+p^3+...)p = p^2+2p^3.
$$ 
The path from the definition of $B_\lfloor^{v,\rightarrow}$
could have length 2 (2 choices) or length 3 (2 choices) or length $\ge 4$. 
In the last case we have $2 \cdot 2 \cdot 2 \cdot 3$ choices for the first 
4 steps, since the first step is left or right, the second step 
cannot backtrace and cannot go to $V_0$ or $V_3$, similarly for the third step, and the fourth step cannot backtrace. Thus 
$$
\P_p(B_\lfloor^{v,\rightarrow})
\le 2p^2+2p^3+ 24 p^4. 
$$
This gives 
\begin{align*}
\P_p(B_\lfloor) &\le  k  \P_p(B_\lfloor^{v,\uparrow} \circ B_\lfloor^{v,\rightarrow})
\le  k (p^2+2p^3)(2p^2+2p^3+ 24 p^4) \\
&= kp^4(2+6p+28p^2+48p^3) \le (2+7.6 p)kp^4. 
\end{align*}
Here we have used a union bound, the BK-inequality, the above estimates
and $p \le \frac 1 {20}$. 
By a similar argument 
\begin{align*}
\P_p(B_\parallel) &\le k(p^2+2p^3)p (k-1) (p^2+2p^3)\\
&=k(k-1)p^5(1+4p+4p^2) \le 2 k(k-1)p^5. 
\end{align*}
On $B_\vdash \cap B_\parallel^c$ we have an open tree with three 
leaves $v,v' \in V_1$ and $v'' \in V_3$ such that $e_v$ is open. 
This tree consists of a path in $V_{2,h}$ of some length $l \ge 1$, 
two bonds $vw, v'w' \in V_{2,v}$ and a bond $v''w'' \in V_{3,v}$,
so that $w,w',w''$ are on the path, among them both endpoints of the path. 
In case of $w'' \neq w,w'$ we have $k$ choices to fix the horizontal path of a given length $l \ge 2$, $l-1$ choices for an interior point of the path, $3!$ choices to distribute $w,w',w''$ among the two endpoints and the interior point, and here we may assume that the two bonds in $E_{3,v}$ adjacent to $w,w'$ are closed (since we are in $B_\parallel^c$). 
In case of $w'' = w$  we have $k$ choices to fix the horizontal path 
of a given length $l \ge 1$, we have $2$ choices to distribute $w,w'$, 
and here we may assume that the bond in $E_{3,v}$ adjacent $w'$ is closed (and similarly for $w'' = w'$). See Figure \ref{Fig:tree} for an illustration of these cases. 
We thus obtain 
\begin{align*}
\P_p(B_\vdash \cap B_\parallel^c)
&\le \sum_{l \ge 2} k(l-1)6 p^{l+4}(1-p)^2 + 2 \sum_{l \ge 1} k 2 p^{l+4} (1-p)\\  
&= kp^5(6p+4) \le 4.3 kp^5
\end{align*}
Using \eqref{equ:inclusion} and combining all the above estimates we get 
\begin{align*}
\P_p(B')  &\le (2+7.6p+4.3p+2 (k-1)p)kp^4 \le (2+10p + 2kp)kp^4.
\end{align*}
 \qed

\begin{Prop}
\label{Prop:estsmall}
\eqref{equ:monopatterngen} is satisfied for all 
$$
p \le \frac 1 {k^2(1+\frac {2} {\ln k})} \quad \text{ and } \quad 
n \ge N_0(k) := 3 k^2(1+\frac 2 {\ln k}).
$$
\end{Prop}


\Pf Let $s_k := k^2(1+ \frac 2 {\ln k })$, and let 
$p \le \frac 1 {s_k}$ and $n \ge \lfloor 3s_k - 4 \rfloor$.  
We note that $n \ge 3k^2-5 \ge \frac{k+2} 2 \ge n_k$, 
and by Lemma \ref{Lem:sk} we have 
$kp \le \frac k {s_k} \le \frac 3 {25}$
and $p \le \frac 1 {s_k} \le \frac 1 {25}$, 
 which gives $2 + 10 p+2kp \le  e$. 
 Thus Lemma \ref{Lem:estsmallp} implies that for  all  $y,x \in \tM^*$
\begin{align*}
&\P_p(\cX_{n+1}^y = x) = \P_p(\cX_{n+1}^y = x, A_{n+1}) + \P_p(\cX_{n+1}^y = x, A_{n+1}^c)\\
&\le \P_p(\cX_n^y = x) \big( p \frac{n+1} 3 + (e kp^4)^{\frac{n-1} 3}p ^{- \frac{k^2-k+1} 2} (1-p)^{-\frac{k^2+k} 2} p^{-n}\Big)\\
&= \P_p(\cX_n^y = x)
\Big( p \frac{n+1} 3 + (e k)^{\frac{n-1} 3}p ^{\frac {n-4} 3 -
\frac{k^2-k+1} 2} (1-p)^{-\frac{k^2+k} 2} \Big) , 
\end{align*}
and it suffices to estimate the last bracket by $\le 1$. 
Since $\frac {n-4} 3 \ge s_k - 3 \ge 
\frac{k^2-k+1} 2 $ the bracket is increasing in $p$, 
so it suffices to consider $p = \frac 1 {s_k}$, 
and by Proposition \ref{Prop:monoN} it suffices to consider 
$n =  \lfloor 3 s_k - 4 \rfloor$. Since 
clearly $ekp < 1$ it suffices to show that 
\begin{align*}
&\frac 1 {s_k} \frac{3 s_k - 3} 3 + (e k)^{\frac{3s_k - 6} 3}
s_k ^{\frac{k^2-k+1} 2 - \frac {3s_k - 9} 3} (1-\frac 1 {s_k})^{-\frac{k^2+k} 2} \le 1 \\
&\Lra \quad 
(e k)^{s_k - 2}
s_k ^{\frac{k^2-k+9} 2  - s_k} (1-\frac 1 {s_k})^{-\frac{k^2+k} 2} \le 1 \quad 
\Lra \quad f(k) \ge 0, \text{ where }\\
&f(k):= (s_k - \frac{k^2-k+9} 2) \ln s_k -
(s_k - 2) \ln (ek) + \frac{k^2+k}2 \ln(1 - \frac 1 {s_k}).
\end{align*}
The latter is shown in Lemma \ref{Lem:sk}.  \qed

\subsection{Large values of $p$}
 
Again, the general strategy in this case is similar to the one for general graphs in \cite{KR}: 
If $p$ is large and we want to go from a pattern $y$ to a pattern $x$ in a large number of steps, 
in many layers all bonds will be open. 
In particular, with high probability there will be such a layer, 
and we proceed by discarding the layer before that. 
This will let us compare  $\P_p(\cX^y_{n+1} = x)$ and $\P_p(\cX^y_n = x)$. 
Again, in comparison to \cite{KR} our estimates need to be done more careful.

\begin{Lem}\label{Lem:estlargep}
For all $n \ge n_k$, $p \in (0,1)$ and $y,x \in \tM^*$ 
we have 
\begin{align} 
\label{equ:estlargepfull}
&\P_p(\cX_{n+1}^y = x, A_{n+1}) \le (1-(1-p)^{k}) \P_p(\cX_n^y = x),\\
\label{equ:estlargepnotfull}
&\P_p(A_{n+1}^c) \le (1-p^{2k})^n,
\\
\label{equ:estlargepprobvar}
&\P_p(\cX_n^y = x) \ge p^{\frac{k^2-k+1}2}(1-p)^{\frac{k^2+k}2} p^{n}
,  
\end{align}
where $A_{n+1} := \bigcup_{2 \le i \le n+1} B_i$ and 
$B_i := \{\forall e \in E_i: \cZ_e = 1\}$. 
\end{Lem}

\Pf For \eqref{equ:estlargepfull} let $B_{i+1}^{n+1} := \bigcap_{i<j \le n+1} B_i^c$. This enables us to write $A_{n+1}$ as the disjoint union
$A_{n+1} =   \bigcup_{2 \le i \le n+1} B_i \cap B_{i+1}^{n+1}$
(considering the topmost layer with all edges open). 
Similarly as in the proof of Lemma \ref{Lem:estsmallp} we have 
\begin{align*}
&\P_p(\cX_{n+1}^y = x, A_{n+1}) = 
\sum_{i,y'} \P_p(\cX_{i-2}^y = y', B_i, \cX^y_i = x_*,  B_{i+1}^{n+1},  \cX_{n+1}^y = x) \\
&= \sum_{i,y'} \P_p(\cX_{i-2}^y = y')\P_p(\cX_2^{y'} = x_*,B_2)\P_p(B_{1}^{n+1-i},\cX_{n+1-i}^{x_*} = x), 
\end{align*}
where $2 \le i \le n+1$ and $y' \in \tM^*$, and in the first step 
we have used a union bound noting that 
$\{\cX^y_{n+1} \in \tM^*\} \cap B_i \subset \{\cX_i ^y = x_*\}$, 
and in the second step we have decomposed the event 
into parts depending on disjoint sets of bonds and used independence. 
For the middle term we have   
\begin{align*}
\P_p(\cX_2^{y'} = x_*,B_2) &= 
\P_p(\cX^{y'}_1 \in \tM^*,B_2) = \P_p(\cX^{y'}_1 \in \tM^*) \P_p(B_2)\\
&\le (1-(1-p)^k) \P_p(B_1), 
\end{align*}
where in the second step again we have decomposed the event into parts depending on disjoint sets bonds, and in the last step we have used 
that if an infection is transmitted from one layer to the next, 
at least one vertical bond has to be open. Thus we obtain 
\begin{align*}
&\P_p(\cX_{n+1}^y = x, A_{n+1}) \\
&\le 
 \sum_{i,y'} \P_p(\cX_{i-2}^y = y')(1-(1-p)^k) \P_p(B_1)\P_p(B_{1,n+1-i}',\cX_{n+1-i}^{x_*} = x)\\
&=   (1-(1-p)^k) \sum_{i,y'} \P_p(\cX_{i-2}^y = y',B_{i-1},B_{i,n}',\cX_{n}^{y} = x)\\
&=   (1-(1-p)^k)\P_p(\bigcup_{1 \le i \le n} B_i, \cX_{n}^{y} = x) 
\le  (1-(1-p)^k)\P_p(\cX_{n}^{y} = x) 
\end{align*}
arguing as above, which gives \eqref{equ:estlargepfull}. 
For \eqref{equ:estlargepnotfull} we note that 
$$
\P_p(A_{n+1}^c) = \P_p(\bigcap_{2 \le i \le n+1} B_i^c)
= (1-p^{2k})^{n},  
$$
by definition of $A_{n+1}$ and $B_i$ and the independence of the $B_i$. 
\eqref{equ:estlargepprobvar} is the same as \eqref{equ:estsmallpprob} (and we didn't use $p \le \frac 1 {20}$ in its proof). \qed

\begin{Prop} \label{Prop:estlarge}
\eqref{equ:monopatterngen} is satisfied for all 
$$
p > 1-\frac 1 {2k} \quad \text{ and }    \quad 
n \ge N_4(k,p) := \frac{(k^2+4k) \ln(1-p)}{2\ln(2k(1-p))}.
$$ 
\end{Prop}

\Pf Let $n,p$ be as above. We first note that 
$2k(1-p) < 1$ and $\frac{\ln(1-p)}{\ln(2k(1-p))} > 1$ and thus 
$N_4(k,p) \ge \frac{k^2 + 4k} 2 \ge \frac{k+2} 2 \ge n_k$. 
So for all $y,x \in \tM^*$ the preceding lemma gives 
\begin{align*}
&\P_p(\cX_{n+1}^y = x) \le \P_p(\cX_{n+1}^y = x, A_{n+1}) + \P_p(A_{n+1}^c)\\
&\le \P_p(\cX_n^y = x) \Big((1-(1-p)^k)
+ ((1-p^{2k})^n  p^{-n-\frac{k^2-k+1}2}(1-p)^{-\frac{k^2+k}2} 
\Big),  
\end{align*}
and it suffices to show that the last bracket is $\le 1$, i.e.
\begin{align*}
&(\frac{1-p^{2k}} p)^n  \le (1-p)^{\frac{k^2+3k} 2}  p^{ \frac{k^2-k+1} 2}.
\end{align*}
By Lemma \ref{Lem:k2}  $\frac {1 - p^{2k}} p \le 2k(1-p)$
and $p^{k^2-k+1} \ge (1-p)^k$, so it suffices to check that 
$(2k(1-p))^n \le (1-p)^{\frac{k^2+4k}2}$, which holds by 
definition of $N_4(k,p)$.  
 \qed

\subsection{Uniform estimate in $p$}

Here we prove Theorem \ref{Thm:patterngeneral} by collecting the estimates 
of the previous subsections. We first note that for all $k \ge 3$ we have 
$$
0 < \frac 1 {3 k^2} < 0.1 < 0.315 
< \frac 4 9 < \frac 2 3 < 1 - \frac 1 {k^2} < 1.
$$
This divides $(0,1)$ into 7 subintervals, and for each subinterval we choose an appropriate estimate for $n$ so that \eqref{equ:monopatterngen}
holds as desired, using Proposition~\ref{Prop:estsmall}, 
Proposition~\ref{Prop:estintermediate} (a),(b),(c) or 
Proposition \ref{Prop:estlarge} respectively. It suffices to note that for $k \ge 3$ we have 
\begin{align*}
&\forall p \in (0,\frac 1 {3 k^2}]:  && N_0(k) \le 9 k^2 \le N(k)\\
&\forall p \in [\frac 1 {3k^2}, 0.1]: && N_1(k,p) \le 
 15 k^8 1.53^k \le N(k) \\
&\forall p \in [0.1, 0.315]: && N_1(k,p) 
\le 526 k^{11/2} 1.95^k \le N(k),\\
&\forall p \in [0.315, \frac 4 9]: && N_2(k,p) 
\le  2119 k^{9/2} 1.95^k  \le N(k),\\
&\forall p \in [\frac 4 9, \frac 2 3]: && N_3(k,p) 
\le 209 k^2 1.95^k \le N(k), \\
&\forall p \in [\frac 2 3, 1-\frac 1 {k^2}]: && N_3(k,p) 
\le 16 k^{5/2} 1.5^k \le N(k), \\
&\forall p \in [1-\frac 1 {k^2},1): && N_4(k) \le 7k^2 \le N(k). 
\end{align*}
Here in the first case the two estimates are easily verified
using $1 + \frac 2 {\ln k} \le 3$. 
For the other cases we refer 
to Lemma \ref{Lem:optimization1} - Lemma \ref{Lem:optimization6} respectively. 
\qed

\section{Number of infected points per layer} \label{Sec:Number}

Next we will prove Theorem \ref{Thm:numberZ2}. 
Let $G = \Z$, $n \ge 0$ and $p \le 0.35$. 
The corresponding general result in \cite{KR} makes use of the recursive estimate 
$$
\E_p(\cW_n) \max_{w \in G} \E_p(\tcW_1^w) \ge \E_p(\cW_{n+1}), 
$$ 
where $\tcW_1^w = |\{v \in V: (v,1) \lra_{\ge 1} (w,0)|$ 
and $\lra_{\ge 1}$ is the connectivity relation 
induced by percolation on $E_{[1,\infty)}$ (see (6.1) of \cite{KR}). 
Here we use the same recursive estimate, but improve on the bound of $\E_p(\tcW_1^w)$. 
Using the recursive estimate, symmetry and monotonicity in $p$, 
it suffices to check that   $p \E_p(\tcW_0) \le 1$ for $p = 0.35$, 
where 
$\tcW_0 = |\{v \in V: (v,0) \lra (0,0) \text{ in } E_{(0,\infty)}\}|$. 
Let $W_{l}$ denote the set of (vertex-)self-avoiding paths in $E_{(0,\infty)}$ 
of length $l$ from $(0,0)$ to some vertex in $V_0$. 
For $P \in W_l$ the size of $P$ is the number of squares enclosed 
between $P$ and the straight path in $V_0$ that starts and ends in the same vertices as $P$. 
We note that on $\{(0,0) \lra (n,0)\}$ there is a unique 
smallest  (i.e. with minimal size) open path connecting $(0,0)$ and $(n,0)$. 
(For clarification: This smallest path does not necessarily have minimal length among all 
open connecting paths.) Let $A_P$ denote the event that $P$ is the 
smallest open path connecting the endpoints of $P$. 
For $l \ge 1$ two paths in $W_l$ are contained in $V_0$, and we let $W_l'$ denote the set of remaining paths. 
We note that 
$$
\E_p(\tcW_0) = \sum_{l \ge 0} \sum_{P \in W_l} \P_p(A_P) = 
 1 + \sum_{l \ge 1} 2 p^l +  \sum_{l \ge 3}p_l',  
\text{ where } 
p_l' =  \sum_{P \in W_l'} \P_p(A_P)  
$$
using $W_l' = \emptyset$ for $l \in \{0,1,2\}$. 
We are left to estimate $p_l'$ for $l \ge 3$. For 
$l \in \{3,4,5\}$ we compute  (see Figure \ref{Fig:pathcounting4}) 
\begin{align*}
&p_3'= 2 p^3(1-p), \quad p_4'= 4p^4(1-p) + 2p^4((1-p)^3+3p(1-p)^2) \quad \text{ and }\\
&p_5'= 6p^5(1-p)^2 + 6 p^5(1-p) + 4 p^5 ((1-p)^3+3p(1-p)^2)\\
&\qquad + 2 p^5((1-p)^5+ 5p(1-p)^4+ 8p^2(1-p)^3). 
\end{align*}
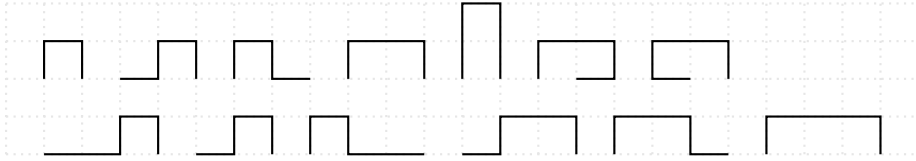
\begin{figure}[htb!] 
\centering
\vspace*{-0.2 cm}
\begin{tikzpicture}[scale = 0.5]
\draw[dotted, thick, lgray] (0,0) grid (24,4);
\draw[thick]
(1,2)--(1,3)--(2,3)--(2,2)
(3,2)--(4,2)--(4,3)--(5,3)--(5,2)
(6,2)--(6,3)--(7,3)--(7,2)--(8,2)
(9,2)--(9,3)--(11,3)--(11,2)
(12,2)--(12,4)--(13,4)--(13,2)
(14,2)--(14,3)--(16,3)--(16,2)--(15,2)
(18,2)--(17,2)--(17,3)--(19,3)--(19,2)
(1,0)--(3,0)--(3,1)--(4,1)--(4,0)
(5,0)--(6,0)--(6,1)--(7,1)--(7,0)
(8,0)--(8,1)--(9,1)--(9,0)--(11,0)
(12,0)--(13,0)--(13,1)--(15,1)--(15,0)
(16,0)--(16,1)--(18,1)--(18,0)--(19,0)
(20,0)--(20,1)--(23,1)--(23,0)
;
\end{tikzpicture}
\caption{Self-avoiding paths of length $3 \le l \le 5$, 
which start in $(0,0)$, end in $V_0$ to the right of $(0,0)$, 
and are contained in $E_{(0,\infty)}$ (but not in $E_0^h$). 
If such a path is  the smallest open path connecting its endpoints, 
all edges of the paths need to be open and some of the edges enclosed by the path need to be closed, 
e.g. for the fourth path exactly two or three of the three enclosed edges need to be closed, 
and for the last path exactly five, four or three of the enclosed edges need to be closed, and 
in the latter case it can't be the three rightmost or the three leftmost edges. 
This gives $p_l'$ for $l \in \{3,4,5\}$. 
}
\label{Fig:pathcounting4}
\end{figure} 

For $l \ge 6$ and every path $P \in W_l'$ we consider its edges in 
the top most layer, which form a set of subpaths. 
Depending on whether the leftmost of these subpaths has length $1,2$ or $\ge 3$, we get   
$\P_p(A_P) \le p^l(1-p)$, $\P_p(A_P) \le p^l((1-p)^3 + 3p(1-p)^2)$ 
or $\P_p(A_P) \le p^l(1-p^2)^2$ respectively. 
Here the terms the brackets arise due to minimal size of the path, 
similarly to the cases considered in Figure \ref{Fig:pathcounting4}.  
Thus 
$$
\P_p(A_P) \le p^l \max\{1-p,(1-p)^3 + 3p(1-p)^2, (1-p^2)^2\}
= p^l (1-p^2)^2, 
$$
where the last step is easy to verify for all $p \in (0,\frac 1 2)$. 
We thus obtain 
$$
p_l' \le (|W_l|-2) p^l (1-p^2)^2 \quad \text{ for all } l \ge 6. 
$$
While for $6 \le l \le 22$ we will use the precise value of $a_l := |W_l|$, 
for $l \ge 23$ we aim at estimating $a_l$.  
Let $b_l$ denote the number of self-avoiding paths in $E_{(0,\infty)}$ of length $l$ from $(0,0)$ to any vertex, 
let $c_l$ denote the number of self-avoiding paths in $E_{(-\infty,\infty)}$ of length $l$ from $(0,0)$ to any vertex, 
for $v \in N := \{(0,1),(0,-1),(1,0),(-1,0)\}$ and $k \in \{-l,...,l\}$ let 
$a_{l,k}^v$ denote the number of such paths that end in $V_k$ and do not 
contain $v$, and let $d_l := \max\{a_{l,k}^v: v \in N, -l \le k \le l\}$.
The sequence $c_l$ can be found in the OESIS (sequence A001411), and its 
values for $l \le 39$ are given in \cite{CEG}.
For small $l$, the values of $a_l,b_l,d_l$ can easily be determined by hand or using a computer. 
The values below have been computed using short recursive python programs that build paths step by step, 
in each step checking whether the path constructed so far has the desired properties, 
and after the last step checking whether 
the last vertex of the path has the desired properties. 
\begin{align*}
&\begin{tabular}{c||c|c|c|c|c|c|c|c|c|c}
$l$ & 1 & 2 & 3 & 4 & 5 & 6 & 7 & 8 & 9 & 10 \\
\hline
$a_l$ & 2 & 2 & 4 & 8 & 20 & 40 & 100 & 216& 548 & 1224 \\
\hline 
$b_l$ & 3 & 7 & 19 & 49 & 131 & 339 & 899 & 2345 & 6199 & 16225\\
\hline 
$d_l$ & 2 & 4 & 8 &  18 & 40& 90& 218& 516& 1250 & 3090 
\end{tabular}\\
&\begin{tabular}{c||c|c|c|c|c|c|c}
$l$   & 11 & 12 & 13 & 14 & 15& 16 &17 \\
\hline 
$a_l$ & 3112 & 7148 &   18228 & 42696 & 109148 & 259520 & 664868\\
\hline 
$b_l$ &  42811& 112285&  296051 &  777411 &  2049025 &  5384855 &  14190509 
\\
\hline
$d_l$ & 7750 & 19506 &  49184 &  124280 &  314822 & 802458& 2054136
\end{tabular}\\
&\begin{tabular}{c||c|c|c|c|c}
$l$  & 18 & 19 & 20 & 21 & 22\\
\hline 
$a_l$ & 1599448 & 4105276&  9969396 & 25630164& 62724196\\
\hline 
$b_l$ &  37313977 & 98324565 &  258654441& 681552747 & \\
\hline
$d_l$ & 5262230 & 13494874 & 34647816& &
\end{tabular}
\end{align*}
In order to be able to treat larger values of $l$ we note that 
\begin{align}
\label{equ:reca}
&\forall l_1 \ge 1, l_2,l_3 \ge 0:  a_{l_1+l_2+l_3} \le b_{l_1}d_{l_2} b_{l_3} 
\text{ and }\\ 
\label{equ:recd}
&\forall l_1 \ge 1,l_2 \ge 0: d_{l_1+l_2} \le \frac 3 4 c_{l_1} d_{l_2}. 
\end{align}
\eqref{equ:reca} can be seen by decomposing a path contributing to $a_{l_1+l_2+l_3}$ into three paths of lengths $l_1,l_2,l_3$ 
in this order and noting the following: The first and third path 
contribute to $b_{l_1}$ and $b_{l_3}$ respectively, and if these two paths are fixed, 
the second part has to connect two vertices of fixed height, 
and it cannot use the penultimate vertex of the first part; 
thus the number of choices for the second part are $\le a_{l_2,k}^v$ 
for some fixed values of $k,v$. 
\eqref{equ:recd} follows similarly: It suffices to estimate $a^v_{l_1+l_2,k}$ 
for fixed values of $v,k$. We decompose a contributing 
path into two parts of lengths $l_1,l_2$ in this order. The first part 
contributes to $c_{l_1}$, but the first edge is forbidden to go in one of the four possible directions (and by symmetry the four possibilities for the first step contribute equally). For fixed first part, we treat the second part as above. This gives \eqref{equ:recd}. 
For $l \in \{23,...,41\}$ we use \eqref{equ:reca} to estimate 
\begin{align*}
a_l \le a_l', \quad \text{ where } \quad 
&a'_{4k} = b_{2(k-5)}b_{2(k-5)} d_{20},  \quad a_{4k+1}' = b_{2(k-5)}b_{2(k-4)} d_{19}, \\ 
&a_{4k+2}' = b_{2(k-5)}b_{2(k-4)} d_{20},  \quad a_{4k+3}' = b_{2(k-4)}b_{2(k-4)} d_{19}. 
\end{align*}
Here we use the values $b_{2m}$ for $1 \le m \le 6$.  
(It can be seen from the table above that using only $b_l$ with even $l$ 
is advantageous.) 
For $l \ge 42$ we may write $l = 21 m + r$ for $r \in \{42,...,62\}$ and $m \in \{0,1,2,....\}$ 
and choose $l_1 \in \{10,...,21\}, l_2 \in \{12,...,21\}$ 
such that $l_1+l_3 = r-20$. 
Using \eqref{equ:reca} and recursively using \eqref{equ:recd} we thus obtain 
$$
a_l \le b_{l_1}b_{l_3}d_{20 + 21m} \le  b_{l_1}b_{l_3}d_{20}  (\frac 3 4 c_{21})^m \le 
b_{10}b_{12}d_{20} \cdot 2.76^{l-42},  
$$
where in the last step we have used $c_{21} = 2408806028$, which implies 
$\frac 3 4 c_{21} \le 2.76^{21}$, 
as well as $\frac{b_{n+1}}{b_n} \le 2.76$ for all $2 \le n \le 20$, 
which implies $b_{l_1} \le b_{10} 2.76^{l_1 - 10}$ and  $b_{l_3} \le b_{12} 2.76^{l_3 - 12}$. 
Collecting all the above cases we can now estimate 
\begin{align*}
\E_p(\tcW_0) \le 1 &+ \sum_{l \ge 1} 2p^l + \sum_{3 \le l \le 5} p_l' 
+ \sum_{6 \le l \le 22} (a_l-2) p^l(1-p^2)^2\\ 
&+ \sum_{23 \le l \le 41} a_l' p^l(1-p^2)^2 + 
\sum_{l \ge 42} b_{10}b_{12}d_{20} 2.76^{l-42} p^l(1-p^2)^2.
\end{align*}
The two infinite sums are geometric series, and thus computing the right hand side 
gives $p \E_p(\tcW_0) \le 1$ for $p = 0.35$ as desired. 
For $G = C_k$ the argument carries over almost without any change. 
Indeed only some of the equalities above now will be inequalities. \qed

\appendix

\section{Appendix}

Here we collect purely analytical estimates used in the preceding chapters. 
Their proofs  are straightforward, but sometimes somewhat tedious.

\begin{Lem} \label{Lem:ac1}
For all $p \in (0,1)$ we have $(p+(1-p)p^3 + 2(1-p)^2p^5))^4 \ge p c_1(p)^3$. 
\end{Lem}

\Pf Dividing by $p^4$ and setting $(1-p)p^2 =: s \in (0,\frac 4 {27}]$ 
we need to show 
that $(1+s+2s^2)^4 \ge (1+s+3s^2+9s^3)^3$ for all $s \in(0,\frac 4 {27}]$. 
This follows from $1+s+3s^2 + 9s^3 \le (1+s+2s^2)(1+2s^2)$ (which 
is equivalent to $1 - 7s + 4s^2 \ge 0$, which holds for $s \in (0,\frac 4 {27}]$)
and $(1+2s^2)^3 \le 1+s+2s^2$ (which is equivalent to $4s+12s^3+8s^5 \le 1$, 
which holds for $s = \frac 4 {27}$ and thus for all $s \in (0,\frac 4 {27}]$). 
\qed

\begin{Lem} \label{Lem:ae}
For all $a > 0$ we have $(1+\frac 1 a)^{a+\frac 1 2} > e$. 
\end{Lem}

\Pf 
For $h(a) := (a+\frac 1 2) \ln(1+ \frac 1 a)$ we have 
$h'(a) = \ln(1 + \frac 1 a) - \frac{a + \frac 1 2}{a(a+1)}$ and 
$h''(a) = \frac{1}{2a^2(a+1)^2} > 0$ for all $a > 0$, i.e. 
$h$ is strictly convex. 
Since $h(a) \to 1$ for $a \to \infty$ we thus  
get $h(a) > 1$ for all $a > 0$, which implies the claim. \qed

\begin{Lem} \label{Lem:sk}
For all $k \ge 3$ we have  $s_k := k^2(1+\frac 2 {\ln k}) \ge \frac{ 25} 3  k \ge 25$, and 
\begin{align*}
&f(k):= (s_k - \frac{k^2-k+9} 2) \ln s_k -
(s_k - 2) \ln (ek) + \frac{k^2+k}2 \ln(1 - \frac 1 {s_k}) \ge 0. 
\end{align*}
\end{Lem}

\Pf For the first assertion we note that  $h(x) := x (1 + \frac 2 {\ln x})$ 
is increasing for $x \ge 3$
(since $h'(x) = 1 + \frac {2(\ln(x)-1)}{(\ln x)^2}  \ge 0$), 
which implies $\frac{s_k}k  \ge h(3) \ge \frac{25}3 \ge \frac{25} k$. 
The second assertion can easily be checked for $k \in \{3,4,5,6,7\}$. 
For $k \ge 8$ we first note that $s_k \ge k^2+k$ (which is equivalent to $2 k \ge \ln k$) 
and that $g(x) = \frac 1 x \ln(1- x) = 
- \sum_{k \ge 0} \frac 1 {k+1} x^k$ is decreasing for $x \in (0,1)$, 
which implies 
\begin{align*}
&\frac{k^2+k}2 \ln(1 - \frac 1 {s_k}) \ge \frac{k^2+k}2 \ln(1 - \frac 1 {k^2+k}) \ge \ln \frac 1 2 \ge - 2.
\end{align*}
Furthermore $s_k \ge k^2 \ge \frac{k^2-k+9}2$ implies
$$
(s_k - \frac{k^2-k+9} 2) \ln s_k \ge (s_k - \frac{k^2-k+9} 2) \ln (k^2). 
$$ 
Thus 
\begin{align*}
&f(k) \ge \Big(k^2(1 + \frac 2 {\ln k}) - \frac{k^2-k+9} 2\Big) 2 \ln k
- k^2(1 + \frac 2 {\ln k}) (\ln k + 1)  + 2 \ln k \\
&= k^2( 1 - \frac 2 {\ln k }) + (k-7) \ln k  \ge 0. 
\end{align*}
\qed 

\begin{Lem} \label{Lem:k2}
For all $k \ge 3$  and $p \ge 1 - \frac 1 {2k}$ we have 
$$
\frac {1-p^{2k}}{p(1-p)} \le 2k \qquad \text{ and } \qquad 
p^{k^2-k+1} \ge (1-p)^k. 
$$
\end{Lem}

\Pf For the first assertion let 
$f(p) := \frac{1-p^{2k}}{p(1-p)} = \frac 1 p + 1 + p + ... + p^{2k-2}$ and note that 
$f'(p) \ge -\frac 1 {p^2} + 1 + 2p \ge - \frac 1 {(\frac 5 6)^2} + 1 + \frac {10} 6 \ge 0$ (since $p \ge 1 - \frac 1 {2k} \ge \frac  5 6$). Thus $f$ is increasing and since $f(1) = 2k$ we are done. 
For the second assertion it suffices to show that  $p^k \ge 1-p$ 
and by monotonicity in $p$ it suffices to consider $p = 1 - \frac 1 {2k}$, i.e. we need to check that $(1-\frac 1 {2k})^k \ge \frac 1 {2k}$. 
Since this is true for $k = 3$, $\frac 1 {2k}$ is decreasing in $k$,  
and $(1-\frac 1 {2k})^k$ is increasing in $k$ 
(noting that $(1-\frac 1 x)^x= e^{-\sum_{k \ge 1} \frac 1 {kx^{k-1}}}$ is increasing in $x$ for $x > 1$), 
we are done. \qed

\bigskip 

The remaining Lemmas aim at estimating the $p$-dependent terms of $N_i(k,p)$ for $i \in \{1,...,4\}$. As preparation we collect several monotonicity properties. 

\begin{Lem} \label{Lem:Nkpmono}
For $p \in (0,1)$ we have:  
\begin{itemize}
\item[(a)] 
$-\ln(p(1-p))$ is $> 0$, decreasing for $p \in (0,\frac 1 2]$, increasing for $p \in [\frac 1 2,1)$. 
\item[(b)] 
$\frac{1-p(1-p)}{1-p}$ is $> 1$, increasing in $p$.
\item[(c)] 
$\frac{c_1(p)} p$ is $> 0$, increasing for $p \in (0,\frac 2 3]$, decreasing for $p \in [\frac 2 3,1)$. 
\item[(d)] $\frac{c_2(p)}p$ and $\frac{c_2(p)-p}{p^2}$ are $> 0$, increasing for $p \in (0,\frac 1 3)$. 
\item[(e)] 
$c_3(p)$ is $> p$, increasing in $p$. 
\end{itemize}
\end{Lem}

\Pf 
(a) follows from the fact that $- \ln(x)$ is  $> 0$ and decreasing for $x \in (0,1)$. 
(b)  follows from $\frac{1 - p(1-p)}{1-p}= 1 + \frac{p^2}{1-p}$. 
(c) follows from $\frac{c_1(p)}{p} = h((1-p)p^2)$ for 
$h(s) = 1+s+3s^2+9s^3$, where $h$ is increasing for $s \in (0,1)$ and 
$(1-p)p^2$ is increasing for $p \in (0,\frac 2 3]$ and decreasing for $p \in [\frac 2 3,1)$. 
(d) follows from $\frac{c_2(p)-p}{p^2} = 2+2p+2p^2+4p^3+8p^4+\frac{64}{1-3p}p^5$ and similarly for $\frac{c_2(p)}p$. 
(e) follows from $c_3(p)  = \Big(\sqrt{\frac 1 {p^2} + \frac {e^2} 4} - \frac e 2\Big)^{-1}$.  
\qed

\begin{Lem} \label{Lem:optimization1}
For $p \in [\frac 1 {3 k^2},0.1]$ and $k \ge 3$ 
we have 
\begin{align*}
N_1(k,p) \le N_{1}(k) := 15 k^8 1.53^k \le N(k). 
\end{align*}
\end{Lem}

\Pf In order to estimate the constituent parts of 
$$
N_1(k,p) = 
 \frac{- \ln(p(1-p))(\frac k {2ec_3(p)} +5)(k^2+ \frac 3 2 k+2) k^{5/2}  (\frac{c_2(p)}{p})^{\frac k {2ec_3(p)} +4}}{
(\frac{1-p}{1-p(1-p)})^{\frac k 2} (1-p) c_1(1-p)^{k-2}}
$$
we use Lemma \ref{Lem:Nkpmono}, $p \in [\frac 1 {3 k^2},0.1]$ and $k \ge 3$. Thus we obtain 
\begin{align*}
&-\ln (p(1-p)) \le -\ln (\frac 1 {3k^2}(1-\frac 1 {3k^2}))
\le 2 \ln k + 1.2 \le 2 \sqrt k, \\
&\frac k {2ec_3(p)} +5 \le \frac k {2ep} + 5 \le 
\frac {3k^3}{2e} + 5 \le k^3(\frac 3 {2e} + \frac 5 {27}) \le 0.74 k^3\\
&(\frac{c_2(p)}{p})^{\frac 1 {2ec_3(p)}} \le (1+2.23p)^{\frac 1 {2ep}}
\le e^{\frac{2.23}{2e}} \le 1.508\end{align*} 
where in the first estimate $\ln k + 0.6 \le \sqrt k$ for $k \ge 4$ 
follows by differentiation. Similarly 
$k^2+ \frac 3 2 k + 2 \le \frac{31}{18}  k^2$, 
$\frac{1-p(1-p)}{1-p} \le 1.012$, $\frac 1 {1-p} \le \frac{10}9$ and 
$\frac 1 {c_1(1-p)} \le 1.006$. 
Using these estimates in the above definition of $N_1(k,p)$ we obtain 
$N_1(k,p) \le N_1(k)$. 
For $N_1(k) \le N(k)$ we use $\frac{1.95^k}{1.53^k} \ge e^{0.242 k}
\ge (0.242 k)^2$ (where the last step follows from 
$\ln x \le \ln 2 + \frac 1 2 (x-2) \le \frac 1 2 x$). 
\qed 

\begin{Lem}\label{Lem:optimization2}
For $p \in [0.1,0.315]$ and $k \ge 3$ 
we have 
\begin{align*}
N_1(k,p) \le N_2(k) := 526 k^{11/2} 1.95^k
\le N(k). 
\end{align*}
\end{Lem}

\Pf 
We estimate $N_1(k,p) \le  kg(p) \frac{31}{18} k^2 k^{5/2} f(p)^k$, where 
\begin{align*}
&f(p) := (\frac{c_2(p)} p)^{\frac 1 {2ec_3(p)}}(\frac {1 - p(1-p)}{1-p})^{\frac 1 2} \frac 1 {c_1(1-p)}\quad \text{ and }\\
&g(p) := - \ln(p(1-p)) (\frac {1} {2ec_3(p)} + \frac 5 3) (1-p)(\frac{c_2(p)}{p})^{4}(\frac{c_1(1-p)}{1-p})^2. 
\end{align*}
For estimating $f$ we divide $[0.1,0.315]$ into the 5 subintervals generated by $\{0.1, 0.18, 0.27, 0.305, 0.313, 0.315\}$, 
and we estimate each of the terms $\frac{c_2(p)} p$, $\frac 1 {2ec_3(p)}$, $\frac {1 - p(1-p)}{1-p}$, $\frac 1 {c_1(1-p)}$ on each subinterval 
using monotonicity in $p$ by Lemma  \ref{Lem:Nkpmono}.  
Thus it is easy to check that $f(p) \le 1.95$. 
Similarly we treat $g(p)$, 
using the 4 subintervals generated by 
$\{0.1,0.25, 0.31,$ $0.3145, 0.315\}$ and obtain  $g(p) \le 305$. 
This gives $N_1(k,p) \le N_1(k)$.  
$N_1(k) \le N(k)$ follows from $\sqrt k \ge \sqrt 3$. 
\qed 

\begin{Lem}\label{Lem:optimization3}
For $p \in  [0.315,\frac 4 9]$ and $k \ge 3$ 
we have 
\begin{align*}
&N_2(k,p) \le N_3(k) :=  2119 k^{9/2} 1.95^k \le N(k). 
\end{align*}
\end{Lem}

\Pf We estimate $N_2(k,p) \le  k g(p) \frac{31}{18} k^2 k^{3/2} f(p)^k$, 
where  
\begin{align*}
&f(p) := \Big( \frac {1 - p(1-p)}{p(1-p)}\Big)^{\frac 1 {2ec_3(p)}} \frac 1 {c_1(1-p)} \quad \text{ and }\\
&g(p) := -\ln(p(1-p)) (\frac 1 {2ec_3(p)} + \frac 5 3 ) \frac{1}{p^5}(1-p(1-p))(\frac{c_1(1-p)}{1-p})^2.
\end{align*}
For estimating $f$ we divide $[0.315,\frac 4 9]$ into the 8 subintervals generated by $\{0.315,$ $0.32,  0.33,0.34,0.36,0.38,0.41,0.43, \frac 4 9\}$, and estimate each of the terms of $f$ on each subinterval using the 
monotonicity relations from Lemma \ref{Lem:Nkpmono}. 
Thus we obtain $f(p) \le 1.95$. 
Similarly we treat $g(p)$ on the interval $[0.315,\frac 4 9]$
and obtain $g(p) \le 1230$. 
Thus we obtain $N_2(k,p) \le  N_3(k)$. 
$N_3(k) \le N(k)$ follows from $k^{\frac 3 2 } \ge 3^{\frac 3 2}$. 
\qed 

\begin{Lem}\label{Lem:optimization4}
For $p \in  [\frac 4 9, \frac 2 3]$ and $k \ge 3$ 
we have 
\begin{align*}
&N_3(k,p) \le  N_4(k) := 209 k^2 1.95^k \le N(k). 
\end{align*}
\end{Lem}

\Pf We note that  
$$
N_3(k,p) =   \frac {-  \ln(p(1-p)) \frac 3 2 (k^2+3k+3)}{p^4 c_1(p)^{k-2}} \le N_4(k) 
$$
again follows from using the monotonicity relations from Lemma \ref{Lem:Nkpmono} on the interval $[\frac 4 9, \frac 2 3]$. 
$N_4(k) \le N(k)$ is obvious. \qed

\begin{Lem}\label{Lem:optimization5}
For $p \in  [\frac 2 3, 1 - \frac 1 {k^2}]$ and $k \ge 3$ 
we have 
\begin{align*}
&N_3(k,p) \le N_5(k) := 16 k^{5/2} 1.5^k \le N(k). 
\end{align*}
\end{Lem}

\Pf Again $N_3(k,p) \le N_5(k)$ follows from
estimating $-\ln(p(1-p)) \le 2 \ln k + \frac 1 {8} \le 2 \sqrt k$, 
$k^2+3k+3 \le \frac 7 3k^2$
and $c_1(p) \ge p \ge \frac 2 3$. 
$N_5(k)\le N(k)$ can easily be seen. \qed

\begin{Lem}\label{Lem:optimization6}
For $p \in  [1 - \frac 1 {k^2},1)$ and $k \ge 3$ 
we have 
\begin{align*}
&N_4(k,p) \le N_6(k) :=  7 k^2 \le N(k). 
\end{align*}
\end{Lem}

\Pf $N_4(k,p) \le N_5(k)$ follows from $\frac{k^2+4k}2 \le \frac 7 6 k^2$ 
and 
$$
\frac{\ln(1-p)}{\ln(2k(1-p))} \le \frac{\ln(k^2)}{\ln(\frac k 2)}
= \frac {2 \ln k}{\ln k - \ln 2} \le \frac{2 \ln 3}{\ln 3 - \ln 2} \le 6
, 
$$
where we have used monotonicity in $p$ and monotonicity in $k$.  
$N_6(k)\le N(k)$ can easily be seen. \qed

\end{document}